\UseRawInputEncoding
\documentclass[11pt,a4paper]{article}
\usepackage{amsfonts}
\usepackage{amssymb}
\usepackage{mathrsfs}
\usepackage{amsmath}
\usepackage{booktabs}
\usepackage{epsf,epsfig,amsfonts,amsgen,indentfirst}
\usepackage{amsmath,amstext,amsbsy,amsopn,amsthm,bbding,wasysym}
\usepackage{multicol,mathdots}
\usepackage{subfigure}
\usepackage{appendix}
\allowdisplaybreaks

\newtheorem{theorem}{Theorem}[section]

\newtheorem{lemma}[theorem]{Lemma}

\newtheorem{claim}{Claim}
\newtheorem{Lemma}{Lemma}

\newtheorem{Problem}{Problem}
\newtheorem{Theorem}{Theorem}
\newtheorem{conjecture}[theorem]{Conjecture}
\newtheorem{problem}[theorem]{Problem}

\renewcommand{\thetheoremA}{\Alph{theoremA}}
\renewcommand{\thelemmaA}{\Alph{lemmaA}}
\renewcommand{\theclaimA}{\Alph{claimA}}
\renewcommand{\thefactA}{\Alph{factA}}

\renewcommand{\baselinestretch}{1.35}

\begin{document}
\title
{\LARGE \textbf{Brualdi-Goldwasser-Michael problem for maximum permanents of {\rm(0,1)}-matrices\thanks{Supported by the National Natural Science Foundation of China (No. 12261071).} }}
\author{ Tingzeng Wu$^{a,b}$\thanks{{Corresponding author.\newline
\emph{E-mail address}: mathtzwu@163.com, a3566293588@163.com, lvhz01000000uestc.edu.cn}}, Xiangshuai Dong$^a$, Huazhong L$\ddot{u}^c$ \\
{\small $^{a}$ School of Mathematics and Statistics, Qinghai Nationalities University, }\\
{\small  Xining, Qinghai 810007, P.R.~China} \\
{\small $^{b}$ Qinghai Institute of Applied Mathematics,   Xining, Qinghai, 810007, P.R.~China}\\
{\small $^{a}$ School of Mathematical Science,  University of Electronic Science and Technology of China,}\\
 {\small Chengdu, Sichuan, 611731, P.R.~China}}
\date{}

\maketitle
\noindent {\bf Abstract:}\ \ Let $\mathscr{U}(n,\tau)$ be the set of all  {\rm(0,1)}-matrices of order $n$ with exactly $\tau$ 0's.  Brualdi et al. investigated the maximum permanents of all matrices in $\mathscr{U}(n,\tau)$(R.A. Brualdi, J.L. Goldwasser, T.S. Michael, Maximum permanents of matrices of zeros and ones, J. Combin. Theory Ser. A 47 (1988) 207--245.). And they
put forward an open problem to characterize the maximum permanents among all matrices in $\mathscr{U}(n,\tau)$. In this paper,  we focus on the problem. And we characterize the maximum permanents of all matrices in $\mathscr{U}(n,\tau)$ when  $n^{2}-3n\leq\tau\leq n^{2}-2n-1$. Furthermore, we also prove  the maximum permanents of all matrices in $\mathscr{U}(n,\tau)$   when $\sigma-kn\equiv0 (mod~k+1)$ and $(k+1)n-\sigma\equiv0(mod~k)$, where $\sigma=n^{2}-\tau$,  $kn\leq\sigma\leq (k+1)n$ and $k$ is integer.

\section{Introduction}%
Let $M=[a_{ij}]_{n\times n}$ be a (0,1)-matrix. The permanent of $M$ is defined by
$$ {\rm per}M= \sum\limits_{\phi}\prod\limits_{i=1}\limits^{n} a_{i}\phi(i),$$
where the summation takes all permutation $\phi$ of $\{1,2,3,\cdots,n\}$. Valiant \cite{val} has shown that computing the permanent is $\#P$-complete even  when restricted to {\rm(0,1)}-matrices.

The maximum permanent of all {\rm(0,1)}-matrices was first considered by  Ryser \cite{rys1}.  And  he put forward a conjecture  as follows.

\begin{conjecture}\label{conj1.1}(Ryser Conjecture, \cite{rys1})
Let $\mathscr{U}(v,k)$ denote all {\rm(0,1)}-matrices of order $v$ with exactly $k$
1's in each row and column. Let $J$ denote the all-$1$ matrix of order $k$. Let $k$ divide $v$ and let $J^{*}$ denote the direct sum of $v/k$ matrices $J^{*}$. Then $J^{*}$ belongs to $\mathscr{U}(v,k)$ and the permanent of $J^{*}$ is maximal in $\mathscr{U}(v,k)$.
\end{conjecture}
Minc \cite{min} studied the lower and upper bounds for  permanents of {\rm(0,1)}-matrices. And he put forward   a conjecture which is a generalization of Ryser conjecture as follows.

\begin{conjecture} (\cite{min})\label{conj1.2}
If $A=[a_{ij}]$ is a square {\rm(0,1)}-matrix of order $n$. Then
 \begin{eqnarray*}
 \rm perA\leq\prod\limits_{i=1}\limits^{n}(r_{i}!)^{\frac{1}{r_{i}}},
 \end{eqnarray*}
equality holds if and only if there exist permutation matrices $P$ and $Q$
such that $PAQ$ is a direct sum of matrices all of whose entries are $1$.
\end{conjecture}
Br\`{e}gman\cite{bre} proved Conjecture \ref{conj1.2}, it implies that  Conjecture \ref{conj1.1} is solved.
Based on Conjectures \ref{conj1.1} and  \ref{conj1.2}, Brualdi, Goldwasser and Michael \cite{bru2} considered also the maximum permanents of all  {\rm(0,1)}-matrices. And  they investigated a less restrictive problem for Ryser's Conjecture, i.e.,
\begin{problem}\label{prob1} (\cite{bru2})
 The average number of $1$'s per row and per column of the matrix $A$ is $k$.  Determine the maximum value of permanent of the matrix $A$.
\end{problem}
For convenience, the notations and definitions of this paper are consistent with  the reference \cite{bru2}. The problem \ref{prob1} has not been solved so far.
 Let $\mu(n,\tau)=$max$~\{{\rm per}A:A\in\mathscr{U}(n,\tau)\}$, and let $\sigma$~($=n^{2}-\tau$) be the number of $1$'s in the {\rm(0,1)}-matrix.  Brualdi et al. determined the maximum permanents of all matrices in $\mathscr{U}(n,\tau)$  when $0\leq\tau\leq 2n$ and $n^{2}-2n\leq\tau\leq n^{2}-n$ as follows.

\begin{theorem} (\cite{bru2}) Let $n$ and $\tau$ be integers with $n\geq3$ and $n^{2}-2n\leq\tau\leq n^{2}-n$. Then $u(n,\tau)= 2^{\lfloor\frac{\sigma-n}{2} \rfloor}$. Moreover, for $A\in\mathscr{U}(n,\tau)$, we have ${\rm per}A=\mu(n,\tau)$ if and
only if $A$ is combinatorially equivalent to one of the following matrices in $\mathscr{U}(n,\tau)$:

$F_{n,\sigma}=J_{2}\oplus J_{2}\oplus\cdots\oplus J_{2}\oplus I_{e}$   ~~~~if $\sigma-n$ is even;

$F_{n,\sigma}^{*}$: which is obtained from $F_{n,\sigma}$ by changing an arbitrary $0$ to a $1$   ~~~~if $\sigma-n$ is odd;

$U_{n,\sigma}=(J_{3}-I_{3})\oplus J_{2}\oplus\cdots\oplus J_{2}\oplus I_{e-1}$   ~~~~ if $\sigma-n$ is odd;

$V_{n,\sigma}=(\left[
\begin{matrix}
0 & 1 & 1  \\
1 & 1 & 1  \\
1 & 1 & 1  \\
\end{matrix}
\right]\oplus I_{2})\oplus\cdots\oplus J_{2}\oplus I_{e+1} $   ~~~~if $\sigma-n$ is odd,

where $e=n-2{\lfloor\frac{\sigma-n}{2} \rfloor}$.
\end{theorem}

\begin{theorem}(\cite{bru2})
 Let $n$ and $\tau$ be integers with $n\geq8$ and $0\leq\tau\leq 2n$. Let $A$ be
a matrix in $\mathscr{U}(n,\tau)$ such that ${\rm per}A=\mu(n,\tau)$. Then $A$ is combinatorially equivalent to one of the matrices below:

$(J_{\tau}-I_{\tau})\oplus_{c} J_{n-\tau}$  ~~if $0<\tau\leq n$,

$( I_{1}\oplus J_{2})^{\langle\tau-n\rangle}\oplus_{c} (J_{4n-3\tau}-I_{4n-3\tau})$    ~~if $n<\tau\leq \frac{4n}{3}$,

$( P_{3})^{\langle3\tau-4n\rangle}\oplus_{c} (I_{1}\oplus J_{2})^{\langle3n-2\tau\rangle}$   ~~if $\frac{4n}{3}<\tau\leq \frac{3n}{2}$,

$( O_{2})^{\langle\tau-\frac{3n}{2}\rangle}\oplus_{c}( P_{3})^{\langle2n-\tau\rangle} $  ~~if $ \frac{3n}{2}<\tau\leq 2n$ and $n$~ is even,

$( O_{2})^{\langle\tau-\frac{3n+1}{2}\rangle}\oplus_{c}( P_{3})^{\langle2n-2-\tau\rangle} \oplus_{c}P_{4}\oplus_{c}P_{4} ^{T}$  ~~if $ \frac{3n+1}{2}<\tau\leq 2n-2$ and $n$ is odd,

$( O_{2})^{\langle\tau-\frac{n-3}{2}\rangle}\oplus_{c} P_{5}$ ~~if $\tau=2n-1$ and $n$ is odd,

$( O_{2})^{\langle\tau-\frac{n-3}{2}\rangle}\oplus_{c} I_{3}$  ~~if $\tau=2n$ and $n$ is odd.
\end{theorem}

Furthermore, Brualdi \cite{bru1} also discussed matrix forms under some special conditions as follows.

\begin{theorem}(\cite{bru1}) For integers $k$ and $n$ integers with $0 \leq k \leq n$, we use the more concise notation $\mathscr{A}(n,k)$ for the class of all {\rm(0,1)}-matrices of order $n$ with exactly $k$ 1's in each row and column.

(i) A matrix $A$ in $\mathscr{A}(n,2)$ achieves the maximal permanent if and only if its rows and columns can be permuted to yield $J_{2}\oplus J_{2}\oplus \cdots\oplus J_{2}$~~($n$ is even), and $(J_{3}-I_{3})\oplus J_{2}\oplus \cdots\oplus J_{2}$~~($n$ is odd).

(ii) If $n\equiv1(mod~3)$, then A matrix $A$ in $\mathscr{A}(n,3)$ achieves the maximal permanent if and only if its rows and columns can be permuted to yield $(J_{4}-I_{4})\oplus J_{3}\oplus \cdots\oplus J_{3}$.
\end{theorem}

\begin{theorem} (\cite{bru1}) The matrices in $\mathscr{A}(n,n-2)$achieving the maximal permanent are:

$n\leq4$, $J_{n}-(I_{n}+C_{n})$;

$n=2t+1\geq4$, $J_{n}-(I_{t}+C_{t})\oplus(I_{t+1}+C_{t+1})$.
\end{theorem}

In this paper, we focus on the maximum permanents of all matrices in $\mathscr{U}(n,\tau)$. And we will prove the following main theorems in
Sections 2 and 3.
\begin{theorem}\label{the1}
Let $n$ and $\tau$ be integers with $n\geq3$. Then
$$\mu(n,\tau)=((k+1)!)^{\frac{\sigma-kn}{k+1}}\times(k!)^{\frac{(k+1)n-\sigma}{k}}.$$
Moreover, for $A\in\mathscr{U}(n,\tau)$,  ${\rm per}A=\mu(n,\tau)$ if and only if $A$ is combinatorially equivalent to the following matrix $K_{n,\sigma}$ in $\mathscr{U}(n,\tau)$, where
$$K_{n,\sigma}=J_{k+1}\oplus J_{k+1}\oplus\cdots\oplus J_{k+1}\oplus J_{k}\oplus \cdots\oplus J_{k},$$
 $kn\leq\sigma\leq (k+1)n$, $k$ is an integer, $\sigma-kn\equiv0(mod~k+1)$ and $(k+1)n-\sigma\equiv0(mod~k)$.
\end{theorem}

\begin{theorem}\label{the2}
Let $n$ and $\tau$ be integers with $n\geq8$ and $n^{2}-3n\leq\tau\leq n^{2}-2n-1$ such that $2n+1\leq\sigma\leq 3n$. Then $\mu(n,\tau)=k\cdot6^{\lfloor\frac{\sigma-2n}{3}\rfloor}\times2^{\lfloor\frac{3n-\sigma}{2}\rfloor}$, where $k\in\{1, \frac{3}{2}, 2, 3, 4\}$.

$({\rm\uppercase\expandafter{\romannumeral1}})$ Set $\sigma-2n\equiv0(mod~3)$.

(i) If $3n-\sigma$ is even, then
$$\mu(n,\tau)={\rm per}F_{n,\sigma}= 6^{\lfloor\frac{\sigma-2n}{3} \rfloor} \times2^{\lfloor\frac{3n-\sigma}{2} \rfloor},$$
and $F_{n,\sigma}$ is combinatorially equivalent to  the matrix below:
$$J_{3}\oplus J_{3}\oplus\cdots\oplus J_{3}\oplus J_{2}\oplus \cdots\oplus J_{2}.$$
(ii) If $3n-\sigma$ is odd, then
$$ \mu(n,\tau)={\rm per}A=6^{\lfloor\frac{\sigma-2n}{3} \rfloor} \times2^{\lfloor\frac{3n-\sigma}{2} \rfloor},$$
and $A$ is combinatorially equivalent to one of the six matrices: $U_{n,\sigma}$, $V_{n,\sigma}$, $W_{n,\sigma}$, $X_{n,\sigma}$, $Y_{n,\sigma}$ and $F_{n,\sigma}^{*}$, where
$$U_{n,\sigma}=(J_{3}-I_{3})\oplus J_{3}\oplus\cdots\oplus J_{3}\oplus J_{2}\oplus \cdots\oplus J_{2},$$
$$V_{n,\sigma}=(\left[
\begin{matrix}
0 & 1 & 1  \\
1 & 1 & 1  \\
1 & 1 & 1  \\
\end{matrix}
\right]\oplus I_{2})\oplus J_{3}\oplus\cdots\oplus J_{3}\oplus J_{2}\oplus \cdots\oplus J_{2},$$
$$W_{n,\sigma}=\left[
\begin{matrix}
0 & 0 & 1 & 1  \\
1 & 1 & 0 & 1  \\
1 & 1 & 1 & 0  \\
1 & 1 & 0 & 1  \\
\end{matrix}
\right]\oplus J_{3}\oplus\cdots\oplus J_{3}\oplus J_{2}\oplus \cdots\oplus J_{2},$$
$$X_{n,\sigma}=\left[
\begin{matrix}
0 & 0 & 1 & 1  \\
1 & 0 & 1 & 1  \\
1 & 1 & 0 & 1  \\
1 & 1 & 1 & 0 \\
\end{matrix}
\right]\oplus J_{3}\oplus\cdots\oplus J_{3}\oplus J_{2}\oplus \cdots\oplus J_{2},$$
$$Y_{n,\sigma}=\left[
\begin{matrix}
0 & 0 & 1 & 1  \\
1 & 0 & 1 & 1  \\
1 & 1 & 0 & 0  \\
1 & 1 & 1 & 1 \\
\end{matrix}
\right]\oplus J_{3}\oplus\cdots\oplus J_{3}\oplus J_{2}\oplus \cdots\oplus J_{2},$$

$F_{n,\sigma}^{*}$ is obtained from $J_{3}\oplus J_{3}\oplus\cdots\oplus J_{3}\oplus J_{2}\oplus \cdots\oplus J_{2}\oplus I_{1}$ by changing an arbitrary entry of $0$ to an entry of $1$.

$({\rm\uppercase\expandafter{\romannumeral2}})$ Set $\sigma-2n\equiv1(mod~3)$.

(i) If $3n-\sigma$ is even, then
$$\mu(n,\tau)={\rm per}A=\frac{3}{2}\times6^{\lfloor\frac{\sigma-2n}{3} \rfloor} \times2^{\lfloor\frac{3n-\sigma}{2} \rfloor},$$

and $A$ is combinatorially equivalent to one of the matrices: $R_{n,\sigma}$ and $S_{n,\sigma}$, where
$$R_{n,\sigma}=(J_{4}-I_{4})\oplus J_{3}\oplus\cdots\oplus J_{3}\oplus J_{2}\oplus \cdots\oplus J_{2},$$
$$S_{n,\sigma}=(J_{3}\oplus I_{2})\oplus J_{3}\oplus\cdots\oplus J_{3}\oplus J_{2}\oplus \cdots\oplus J_{2}.$$

(ii) If $3n-\sigma$ is odd, then
$$\mu(n,\tau)={\rm per}A=2\times6^{\lfloor\frac{\sigma-2n}{3} \rfloor} \times2^{\lfloor\frac{3n-\sigma}{2} \rfloor},$$
and $A$ is combinatorially equivalent to one of the matrices: $T_{n,\sigma}$ and $P_{n,\sigma}^{*}$, where
$$T_{n,\sigma}=(\left[
\begin{matrix}
0 & 1 & 1  \\
1 & 1 & 1  \\
1 & 1 & 1  \\
\end{matrix}
\right]\oplus I_{1})\oplus J_{3}\oplus\cdots\oplus J_{3}\oplus J_{2}\oplus \cdots\oplus J_{2},$$~
$$Q_{n,\sigma}=(\left[
\begin{matrix}
0 & 1 & 1 & 1 \\
1 & 1 & 1 & 1 \\
1 & 1 & 1 & 1 \\
1 & 1 & 1 & 1 \\
\end{matrix}
\right]\oplus J_{2}) \oplus J_{3}\oplus\cdots\oplus J_{3}\oplus J_{2}\oplus \cdots\oplus J_{2},$$

$P_{n,\sigma}^{*}$ is obtained from $J_{3}\oplus J_{3}\oplus\cdots\oplus J_{3}\oplus J_{2}\oplus \cdots\oplus J_{2}\oplus J_{2}$~ by changing an arbitrary $0$ to a $1$.

$({\rm\uppercase\expandafter{\romannumeral3}})$ Set $\sigma-2n\equiv2(mod~3)$.

(i) If $3n-\sigma$ is even, then
$$ \mu(n,\tau)={\rm per}M_{n,\sigma}=3\times6^{\lfloor\frac{\sigma-2n}{3} \rfloor} \times2^{\lfloor\frac{3n-\sigma}{2} \rfloor},$$

and $M_{n,\sigma}$ is combinatorially equivalent to  the matrix below:
$$(J_{3}\oplus I_{1})\oplus J_{3}\oplus\cdots\oplus J_{3}\oplus J_{2}\oplus \cdots\oplus J_{2}.$$

(ii) If $3n-\sigma$ is odd, then
$$\mu(n,\tau)={\rm per}N_{n,\sigma}=4\times6^{\lfloor\frac{\sigma-2n}{3} \rfloor} \times2^{\lfloor\frac{3n-\sigma}{2} \rfloor},$$

and $N_{n,\sigma}$ is combinatorially equivalent to  the matrix below:
$$\left[
\begin{matrix}
0 & 1 & 1  \\
1 & 1 & 1  \\
1 & 1 & 1  \\
\end{matrix}
\right]\oplus J_{3}\oplus\cdots\oplus J_{3}\oplus J_{2}\oplus \cdots\oplus J_{2}.$$

\end{theorem}

The outline of this paper is shown as follows. In Section 2, we gave the proof of Theorem \ref{the1}. And  we gave the proof of Theorem \ref{the2} in Section 3. In final section, we summarized the main conclusions of this paper.  And we analyzed Brualdi-Goldwasser-Michael problem.

 \section{ Proof of Theorem \ref{the1}}

Before we give the proof of Theorem \ref{the1}, we first present a lemma as follows.

\begin{lemma}( \cite{bru1}) \label{lem1}
Let $m$ and $t$ be integers with $m\geq2$ and $t\geq1$. Then
$$((m+t-1)!)^{\frac{1}{m+t-1}}\times(m!)^{\frac{1}{m}}>((m+t)!)^{\frac{1}{m+t}}\times((m-1)!)^{\frac{1}{m-1}}.$$
\end{lemma}

\textbf{Proof of Theorem \ref{the1}:}  By Lemma \ref{lem1}, we obtain that $\prod\limits_{i=1}\limits^{n}(r_{i}!)^{\frac{1}{r_{i}}}$ ($r_{i}$ is a positive integer) attains its maximum value when $r_{i}=\lfloor\frac{\sigma}{n}\rfloor$ or $\lfloor\frac{\sigma}{n}\rfloor+1$.  Assume that the number of $r_{i}=\lfloor\frac{\sigma}{n}\rfloor$  is $x$, and  the number of $r_{i}=\lfloor\frac{\sigma}{n}\rfloor+1$ is $y$. For any matrix $A\in\mathscr{U}(n,\tau)$, then
\begin{equation}\label{equ1}
\begin{cases}
x+y=n\\

xr+y(r+1)=\sigma.
\end{cases}
\end{equation}
Solving equation (\ref{equ1}), we obtain that $x=nr+n-\sigma$ and $y=\sigma-nr$. This means that $\mu(n,\tau)=max~\{{\rm per}A:A\in\mathscr{U}(n,\tau)\} \leq max \{\prod\limits_{i=1}\limits^{n}(r_{i})!^{\frac{1}{r_{i}}}$: $r_{i}$ is a positive integer, $\sum\limits_{i=1}\limits^{n}r_{i}=\sigma\}=(r_{i}!)^{\frac{x}{r}}\times((r_{i}+1)!)^{\frac{y}{r+1}}$. Since $kn\leq\sigma\leq(k+1)n$, $\sigma-kn\equiv0(mod~k+1)$ and $kn-\sigma\equiv0(mod~k)$, we have
\begin{eqnarray}\label{equ2}
{\rm per}A=\mu(n,\tau)=((k+1)!)^{\frac{\sigma-kn}{k+1}}\times(k!)^{\frac{(k+1)n-\sigma}{k}},
\end{eqnarray}
where $r=k$, $x=nk+n-\sigma$ and $y=\sigma-nk$. Checking (\ref{equ2}), we know that the form of $A$ is  combinatorially equivalent to  the matrix
$J_{k+1}\oplus J_{k+1}\oplus\cdots\oplus J_{k+1}\oplus J_{k}\oplus \cdots\oplus J_{k}.$\qquad \qquad\qquad\qquad\qquad\qquad\qquad\quad\qquad \qquad\qquad\qquad\qquad\qquad\qquad\quad $\Box$

\section{ Proof of Theorem \ref{the2}}

In this scetion, we will give the proof of Theorem \ref{the2}. First, we  introduce the proving idea of Theorem \ref{the2}.

\textbf{Proving idea:} We prove the Theorem \ref{the2} by induction on $n$. The statement is true for $n=8$. Suppose that the conclusion in Theorem \ref{the2} holds  for $n=m(m>8)$.  The following we prove the conclusion in Theorem \ref{the2} holds for $n=m+1$. Since $A\in\mathscr{U}(n,\tau)$ and $2n+1\leq\sigma\leq 3n$,
we consider the three  cases, i.e., $\sigma<3m+2$, $\sigma=3m+2$ and $\sigma=3m+3$.   Firstly, suppose that $\sigma<3m+2$. This implies that some line of $(m+1)$-square matrix $A$ has a single $1$ or two $1$'s. So, $(m+1)$-square matrix $A$ has two forms as follows.
$$A=B=\left[
\begin{matrix}
\begin{array}{c|c}

C        & B_{1}   \\
\hline
0\cdots0 & 1   \\

\end{array}
\end{matrix}\right]
{\rm or}~A=D=\left[
\begin{matrix}
\begin{array}{c|c|c}

1  &  1  &  0\cdots0 \\
\hline
D_{1}  &  D_{2}  &    C     \\
\end{array}
\end{matrix}\right].$$
According to the numbers of $1$'s in the $B_{1}$,  $D_{1}$ as above, we discuss  the forms of $m$-square submatrix $C$ in $(m+1)$-square matrix $B$ and $m$-square submatrix $D(1,1)$, where $D(1,1)$ is a  $m$-square submatrix obtained by deleting first line and first column in $(m+1)$-square matrix $D$. Thereout the form of $(m+1)$-square matrix $A$ is determined. Secondly, assume that $\sigma=3m+2$.  This indicates that some line of $(m+1)$-square matrix $A$ has a single $1$, two $1$'s or three $1$'s.
If some line of $(m+1)$-square matrix $A$ has two $1$'s, then $A=D=\left[
\begin{matrix}
\begin{array}{c|c|c}

1  &  1  &  0\cdots0 \\
\hline
D_{1}  &  D_{2}  &    C     \\
\end{array}
\end{matrix}\right]$.
If some line of $(m+1)$-square matrix $A$ has three $1$'s, then
$ A= E=\left[
\begin{matrix}
\begin{array}{c|c|c|c}

1  &  1  &  1 &  0\cdots0 \\
\hline
E_{1}  &  E_{2}  &  E_{3} &  C   \\
\end{array}
\end{matrix}\right]$.
 If some line of $(m+1)$-square matrix $A$ has a single $1$, then there exists other line of $(m+1)$-square matrix $A$ has $q$~$1$'s, where  $q\geq 4$. So, $A=H=\left[\begin{matrix}
\begin{array}{c|c|c|c|c}

1  &  1  & \cdots& 1 &  0\cdots0 \\
\hline
H_{1}  &  H_{2} & \cdots &  H_{q} &  C   \\
\end{array}
\end{matrix}\right].$
Simlar to the first case, according to the numbers of $1$'s in the $H_{1}$, $D_{1}$,  $E_{1}$ as above, we discuss  the forms of  $m$-square submatrices  $H(1,1)$, $D(1,1)$ and $E(1,1)$ , where $m$-square submatrices $H(1,1)$ and $E(1,1)$ denote the submatrices obtained by deleting first line and first column in $(m+1)$-square matrices $H$ and $E$, respectively. Thereout the forms of the $(m+1)$-square  matrix $A$ are characterized. Finally, suppose that $\sigma=3m+3$.   This means that some line of $(m+1)$-square matrix $A$ has a single $1$, two $1$'s or three $1$'s. If some line of $(m+1)$-square matrix $A$ has a single $1$ or two $1$'s, then there exists some line of $(m+1)$-square matrix $A$ has $q$~$1$'s, where $q\geq4$, then
$A=H=\left[\begin{matrix}
\begin{array}{c|c|c|c|c}

1  &  1  & \cdots& 1 &  0\cdots0 \\
\hline
H_{1}  &  H_{2} & \cdots &  H_{q} &  C   \\
\end{array}
\end{matrix}\right].$
If some line of $(m+1)$-square matrix $A$ has  three $1$'s, then
$A=E=\left[
\begin{matrix}
\begin{array}{c|c|c|c}

1  &  1  &  1 &  0\cdots0 \\
\hline
E_{1}  &  E_{2}  &  E_{3} &  C   \\
\end{array}
\end{matrix}\right]$.
Simlar to the first case, according to the numbers of $1$'s in the   $H_{1}$ and $E_{1}$ as above, we discuss  the forms of  $m$-square submatrices $H(1,1)$ and $E(1,1)$.  Thereout the forms of the $(m+1)$-square  matrix $A$ are characterized.
Combining the arguments above, the proof is complete.

Base on the proving idea as above, we give  a detailed proof of Theorem \ref{the2}.

\textbf{Proof of Theorem \ref{the2}:} Employing Maple 18.0, the permanents of all matrices in $\mathscr{U}(n,\tau)$ are computed. And the maximum permanents are enumerate, see Table \ref{tab1}. It can be known that the statement is true for $n=8$.
\begin{table}[htbp]
\begin{center}
\caption{Matrices $\mathscr{U}(8,\tau)$  with maximum permanents}
\label{tab1}
\setlength{\tabcolsep}{16mm}
\begin{tabular}{@{}llllll@{} } \toprule
$\sigma$ &$\tau$ &$\mu(8,\tau)$ &permanent \\ \midrule
  $17$&$47$& $16$ & $P_{8,17}^{*}$, $T_{8,17}$  \\
  $18$&$46$& $18$ & $M_{8,18}$  \\
  $19$&$45$& $24$ & $U_{8,19}$, $U_{8,19}$  \\
  $20$&$44$& $36$ & $R_{8,20}$  \\
  $21$&$43$& $48$ & $N_{8,21}$  \\
  $22$&$42$& $72$ & $F_{8,22}$ \\
  $23$&$41$& $72$ & $P_{8,21}^{*}$, $Q_{8,21}$ \\
  $24$&$40$& $108$ & $M_{8,24}$  \\
\bottomrule
\end{tabular}
\end{center}
\end{table}
Suppose that the conclusion in Theorem \ref{the2} holds  for $n=m(m>8)$. Nextly, we will prove the conclusion in Theorem \ref{the2} holds  for $n=m+1$. Since $(m+1)$-square matrix $A\in\mathscr{U}(n,\tau)$ and $2n+1\leq\sigma\leq 3n$,
we consider the three  cases, i.e., $\sigma<3m+2$, $\sigma=3m+2$ and $\sigma=3m+3$.

\textbf{Case 1}.
Suppose that $\sigma<3m+2$. This implies that some line of $(m+1)$-square matrix $A$ has a single $1$ or two $1$'s. Thus $(m+1)$-square matrix $A$ has two forms:
$A=B=\left[
\begin{matrix}
\begin{array}{c|c}

C        & B_{1}   \\
\hline
0\cdots0 & 1   \\

\end{array}
\end{matrix}\right]
{\rm or}~A=D=\left[
\begin{matrix}
\begin{array}{c|c|c}

1  &  1  &  0\cdots0 \\
\hline
D_{1}  &  D_{2}  &    C     \\
\end{array}
\end{matrix}\right].$

Now we determine  the form of $(m+1)$-square matrix $B$. Assume that the number of $1$'s in  $B_{1}$ is $r$. For any $\sigma\in (2n, 3n]$,
 we consider the  following three subcases.

\textbf{Subcase 1.1}. Let $\sigma-2(m+1)\equiv0(mod~3)$. By simple calculation,   we obtain that the permanents of $\{ F_{(m+1),\sigma}, U_{(m+1),\sigma}, V_{(m+1),\sigma}, W_{(m+1),\sigma}, X_{(m+1),\sigma}, Y_{(m+1),\sigma}, F_{(m+1),\sigma}^{*}\}$ equal to $6^{\lfloor\frac{\sigma-2(m+1)}{3} \rfloor} \times2^{\lfloor\frac{3(m+1)-\sigma}{2} \rfloor}$. This implies that ${\rm per}A=\mu(m+1,\tau)\geq6^{\lfloor\frac{\sigma-2(m+1)}{3} \rfloor} \times2^{\lfloor\frac{3(m+1)-\sigma}{2} \rfloor}$. By induction hypothesis, for any $m$-square submatrix $C$ of $B$, $m$-square matrix $C$ having all enumerating matrix forms in Theorem \ref{the2} have the maximum values on permanent of the matices.
Checking $m$-square matrix $C$, we know that the number of $1$'s in $m$-square matrix $C$ is $\sigma-(r+1)$.  Now, we consider four subsubcases.

\textbf{Subsubcase 1.1.1}. Assume $\sigma-2(m+1)+1-r\equiv0(mod~3)$. By Laplace Expansion Theorem, we get that ${\rm per}B={\rm per}C$.
By induction hypothesis, we know that ${\rm per}C\leq6^{\lfloor\frac{\sigma-2(m+1)+1-r}{3} \rfloor} \times2^{\lfloor\frac{3(m+1)-\sigma+r-2}{2} \rfloor}$, and
\begin{eqnarray}\label{equ3}
6^{\lfloor\frac{\sigma-2(m+1)}{3} \rfloor} \times2^{\lfloor\frac{3(m+1)-\sigma}{2} \rfloor}\leq{\rm per}A={\rm per}B={\rm per}C\leq6^{\lfloor\frac{\sigma-2(m+1)+1-r}{3} \rfloor} \times2^{\lfloor\frac{3(m+1)-\sigma+r-2}{2} \rfloor}.
\end{eqnarray}
Since $\sigma-2(m+1)\equiv0(mod~3)$ and $\sigma-2(m+1)+1-r\equiv0(mod~3)$, we have $1-r=0,-3,-6,\cdots$, and~$r=1,4,7,\cdots$. By (\ref{equ3}), we have
\begin{eqnarray}\label{equ4}
6^{\lfloor\frac{\sigma-2(m+1)}{3} \rfloor} \times2^{\lfloor\frac{3(m+1)-\sigma}{2} \rfloor}&\leq &6^{\lfloor\frac{\sigma-2(m+1)}{3} \rfloor} \times2^{\lfloor\frac{3(m+1)-\sigma-1}{2} \rfloor}\times6^{\lfloor\frac{1-r}{3} \rfloor} \times2^{\lfloor\frac{r-1}{2} \rfloor}\nonumber\\
&=&6^{\lfloor\frac{\sigma-2(m+1)}{3} \rfloor} \times2^{\lfloor\frac{3(m+1)-\sigma-1}{2} \rfloor}\times\frac{2^{\lfloor\frac{r-1}{2} \rfloor}}{6^{\lfloor\frac{r-1}{3} \rfloor}}\nonumber\\
&\leq &6^{\lfloor\frac{\sigma-2(m+1)}{3} \rfloor} \times2^{\lfloor\frac{3(m+1)-\sigma}{2} \rfloor}\times\frac{2^{\frac{r-1}{2}}}{6^{\frac{r-1}{3}}}\nonumber\\
&=&6^{\lfloor\frac{\sigma-2(m+1)}{3} \rfloor} \times2^{\lfloor\frac{3(m+1)-\sigma}{2} \rfloor}\times(\frac{\sqrt{2}}{\sqrt[3]{6}})^{r-1}\\
&\leq&6^{\lfloor\frac{\sigma-2(m+1)}{3} \rfloor} \times2^{\lfloor\frac{3(m+1)-\sigma}{2} \rfloor}.\nonumber
\end{eqnarray}
By (\ref{equ4}), we get that $r=1$. And by (\ref{equ3}) and (\ref{equ4}), we obtain that  $3(m+1)-\sigma$ and $3(m+1)-\sigma+r-2$ are odd and even, respectively. According to the arguments as above, it can be known that $m$-square matrix $C$ is combinatorially equivalent to $F_{m,\sigma}$. Hence, $(m+1)$-square matrix $A$ is combinatorially equivalent to $F_{(m+1),\sigma}^{*}$.

\textbf{Subsubcase 1.1.2}. Assume $\sigma-2(m+1)+1-r\equiv1(mod~3)$. Since $\sigma-2(m+1)\equiv0(mod~3)$ and $\sigma-2(m+1)+1-r\equiv1(mod~3)$, we get that $1-r=1,-2,-5,\cdots$. This means that  $r=0,3,6,\cdots$. In this case, we will consider  $3(m+1)-\sigma+r-2$ which is even or odd, respectively. Suppose that $3m-\sigma+r+1(=3(m+1)-\sigma+r-2)$ is even. By induction hypothesis, we know that $m$-square  submatrix $C$ satisfies item (II) in Theorem \ref{the2}. So, we know that $k=\frac{3}{2}$, and
\begin{eqnarray}\label{equ5}
 6^{\lfloor\frac{\sigma-2(m+1)}{3} \rfloor} \times2^{\lfloor\frac{3(m+1)-\sigma}{2} \rfloor}\leq{\rm per}A={\rm per}B={\rm per}C= \frac{3}{2}\times6^{\lfloor\frac{\sigma-2(m+1)+1-r}{3} \rfloor} \times2^{\lfloor\frac{3(m+1)-\sigma+r-2}{2} \rfloor}.
 \end{eqnarray}
By (\ref{equ5}) and $\sigma-2(m+1)\equiv0(mod~3)$, we get that
\begin{eqnarray}\label{equ6}
{\rm per}C&=& \frac{3}{2}\times6^{\lfloor\frac{\sigma-2(m+1)}{3} \rfloor} \times2^{\lfloor\frac{3(m+1)-\sigma-1}{2} \rfloor}\times6^{\lfloor\frac{1-r}{3} \rfloor} \times2^{\lfloor\frac{r-1}{2} \rfloor}\nonumber\\
&\leq&\frac{3}{2}\times6^{\lfloor\frac{\sigma-2(m+1)}{3} \rfloor} \times2^{\lfloor\frac{3(m+1)-\sigma-1}{2} \rfloor}\times{6^{\frac{1-r}{3} }}\times{2^{\frac{r-1}{2}}}\nonumber\\
&\leq& \frac{3}{2}\times6^{\lfloor\frac{\sigma-2(m+1)}{3} \rfloor} \times2^{\lfloor\frac{3(m+1)-\sigma}{2} \rfloor}\times\frac{2^{\frac{r-1}{2}}}{6^{\frac{r-1}{3}}}\\
&=&\frac{3}{2}\times6^{\lfloor\frac{\sigma-2(m+1)}{3} \rfloor} \times2^{\lfloor\frac{3(m+1)-\sigma}{2} \rfloor}\times(\frac{\sqrt{2}}{\sqrt[3]{6}})^{r-1}.\nonumber
\end{eqnarray}
Set $r=0$. By (\ref{equ6}), we have
${\rm per}C=\frac{3}{2}\times6^{\lfloor\frac{\sigma-2(m+1)+1}{3} \rfloor}\times2^{\lfloor\frac{3(m+1)-\sigma}{2} \rfloor}\times\frac{1}{2}
<6^{\lfloor\frac{\sigma-2(m+1)+1}{3} \rfloor}\times2^{\lfloor\frac{3(m+1)-\sigma}{2} \rfloor}$.
This contradicts to ${\rm per}C\geq 6^{\lfloor\frac{\sigma-2(m+1)}{3} \rfloor} \times2^{\lfloor\frac{3(m+1)-\sigma}{2} \rfloor}$ in (\ref{equ5}).
Set $r>0$. By (\ref{equ6}), we know that If $r$ gets bigger then ${\rm per}C$ gets smaller. This implies that ${\rm per}C< 6^{\lfloor\frac{\sigma-2(m+1)}{3} \rfloor} \times2^{\lfloor\frac{3(m+1)-\sigma}{2} \rfloor}$ when $r>0$. This contradicts to ${\rm per}C\geq 6^{\lfloor\frac{\sigma-2(m+1)}{3} \rfloor} \times2^{\lfloor\frac{3(m+1)-\sigma}{2} \rfloor}$ in (\ref{equ5}).
Suppose that $3(m+1)-\sigma+r-2$ is odd. By induction hypothesis, we know that $m$-square  submatrix $C$ satisfies item (II) in Theorem \ref{the2}, i.e.,  $k=2$ in Theorem \ref{the2}, $m$-square submatrix $C$ is combinatorially equivalent to one of the matrices \{$T_{m,\sigma}$, $Q_{m,\sigma}$, $P_{m,\sigma}^{*}$\}, and
\begin{eqnarray}\label{equ7}
 6^{\lfloor\frac{\sigma-2(m+1)}{3} \rfloor} \times2^{\lfloor\frac{3(m+1)-\sigma}{2} \rfloor}\leq{\rm per}A={\rm per}B={\rm per}C= 2\times6^{\lfloor\frac{\sigma-2(m+1)+1-r}{3} \rfloor} \times2^{\lfloor\frac{3(m+1)-\sigma+r-2}{2} \rfloor}.
 \end{eqnarray}
Since $\sigma-2(m+1)\equiv0(mod~3)$, we have $6^{\lfloor\frac{\sigma-2(m+1)+1}{3} \rfloor}\times2^{\lfloor\frac{3(m+1)-\sigma}{2} \rfloor}=6^{\lfloor\frac{\sigma-2(m+1)}{3} \rfloor}\times2^{\lfloor\frac{3(m+1)-\sigma}{2} \rfloor}$. Set $r=0$. ${\rm per}A={\rm per}B={\rm per}C=2\times6^{\lfloor\frac{\sigma-2(m+1)+1}{3} \rfloor}\times2^{\lfloor\frac{3(m+1)-\sigma-2}{2} \rfloor}
=6^{\lfloor\frac{\sigma-2(m+1)}{3} \rfloor}\times2^{\lfloor\frac{3(m+1)-\sigma}{2} \rfloor}.$  By matrix forms of $m$-square matrix $C$, we know that $(m+1)$-square matrix $A$ is combinatorially equivalent to one of the matrices \{$V_{(m+1),\sigma}$, $F_{(m+1),\sigma}^{*}$\}. Set $r=3$. By (\ref{equ7}), we have ${\rm per}C=2\times6^{\lfloor\frac{\sigma-2(m+1)-2}{3} \rfloor}\times2^{\lfloor\frac{3(m+1)-\sigma+1}{2} \rfloor}
=\frac{1}{3}\times6^{\lfloor\frac{\sigma-2(m+1)}{3} \rfloor}\times2^{\lfloor\frac{3(m+1)-\sigma}{2} \rfloor}
<6^{\lfloor\frac{\sigma-2(m+1)}{3} \rfloor}\times2^{\lfloor\frac{3(m+1)-\sigma}{2} \rfloor}.$
This contradicts to ${\rm per}C\geq 6^{\lfloor\frac{\sigma-2(m+1)}{3} \rfloor} \times2^{\lfloor\frac{3(m+1)-\sigma}{2} \rfloor}$ in (\ref{equ7}).
Set $r>3$. By (\ref{equ6}), we know that if $r$ is enlarged then ${\rm per}C$ is lessened. This indicates that ${\rm per}C< 6^{\lfloor\frac{\sigma-2(m+1)}{3} \rfloor} \times2^{\lfloor\frac{3(m+1)-\sigma}{2} \rfloor}$, a contradiction.

\textbf{Subsubcase 1.1.3}. Assume $\sigma-2(m+1)+1-r\equiv2(mod~3)$. Since $\sigma-2(m+1)\equiv0(mod~3)$ and $\sigma-2(m+1)+1-r\equiv2(mod~3)$, we get that $1-r=2,-1,-4,\cdots$. This means that $r=2,5,8,\cdots$. Suppose  that $3m-(\sigma-r-1)=3(m+1)-\sigma+r-2$ is even. By induction hypothesis, we know that $m$-square submatrix $C$ satisfies item (III) in Theorem \ref{the2}, i.e.,  $k=3$ in Theorem \ref{the2}, and
\begin{eqnarray}\label{equ8}
6^{\lfloor\frac{\sigma-2(m+1)}{3} \rfloor} \times2^{\lfloor\frac{3(m+1)-\sigma}{2} \rfloor}\leq{\rm per}A={\rm per}B={\rm per}C= 3\times6^{\lfloor\frac{\sigma-2(m+1)+1-r}{3} \rfloor} \times2^{\lfloor\frac{3(m+1)-\sigma+r-2}{2} \rfloor}.
\end{eqnarray}
By (\ref{equ8}), (\ref{equ6}) and simple calculating, it is easy to obtain that if $r\geq 2$ then ${\rm per}C< 6^{\lfloor\frac{\sigma-2(m+1)}{3} \rfloor} \times2^{\lfloor\frac{3(m+1)-\sigma}{2} \rfloor}$. This contradicts to
 ${\rm per}C\geq 6^{\lfloor\frac{\sigma-2(m+1)}{3} \rfloor} \times2^{\lfloor\frac{3(m+1)-\sigma}{2} \rfloor}$ in (\ref{equ8}).
Suppose that $3m-(\sigma-r-1)=3(m+1)-\sigma+r-2$ is odd. By induction hypothesis, we know that $m$-square submatrix $C$ satisfies item (III) in Theorem \ref{the2}. So, we know that $k=4$, and
\begin{eqnarray}\label{equ9}
6^{\lfloor\frac{\sigma-2(m+1)}{3} \rfloor} \times2^{\lfloor\frac{3(m+1)-\sigma}{2} \rfloor}\leq{\rm per}A={\rm per}B={\rm per}C= 4\times6^{\lfloor\frac{\sigma-2(m+1)+1-r}{3} \rfloor} \times2^{\lfloor\frac{3(m+1)-\sigma+r-2}{2} \rfloor}.
\end{eqnarray}
Set $r=2$. By (\ref{equ9}), we have ${\rm per}C=4\times6^{\lfloor\frac{\sigma-2(m+1)-1}{3} \rfloor} \times2^{\lfloor\frac{3(m+1)-\sigma}{2} \rfloor}=\frac{2}{3}\times6^{\lfloor\frac{\sigma-2(m+1)}{3} \rfloor} \times2^{\lfloor\frac{3(m+1)-\sigma}{2} \rfloor}
<6^{\lfloor\frac{\sigma-2(m+1)}{3} \rfloor} \times2^{\lfloor\frac{3(m+1)-\sigma}{2} \rfloor}.$
This contradicts to ${\rm per}C\geq 6^{\lfloor\frac{\sigma-2(m+1)}{3} \rfloor} \times2^{\lfloor\frac{3(m+1)-\sigma}{2} \rfloor}$ in (\ref{equ9}).
Set $r>2$, By (\ref{equ6}), we know that if $r$ is enlarged then ${\rm per}C$ is lessened. This indicates ${\rm per}C< 6^{\lfloor\frac{\sigma-2(m+1)}{3} \rfloor} \times2^{\lfloor\frac{3(m+1)-\sigma}{2} \rfloor}$, a contradiction.

\textbf{Subcase 1.2}. Let $\sigma-2(m+1)\equiv1(mod~3)$. By simple calculation, there exists $k_{1}\in\{\frac{3}{2},2\}\subseteq \{1, \frac{3}{2}, 2, 3, 4\}$ such that the permanents of $\{ R_{(m+1),\sigma}$, $S_{(m+1),\sigma}$, $T_{(m+1),\sigma}$, $Q_{(m+1),\sigma}$, $P_{(m+1),\sigma}^{*}\}$ equal to $k_{1}\times6^{\lfloor\frac{\sigma-2(m+1)}{3} \rfloor} \times2^{\lfloor\frac{3(m+1)-\sigma}{2} \rfloor}$. It implies that $u((m+1),\tau)\geq k_{1}\times6^{\lfloor\frac{\sigma-2(m+1)}{3} \rfloor} \times2^{\lfloor\frac{3(m+1)-\sigma}{2} \rfloor}$. By induction hypothesis, for any $m$-square submatrix $C$ of $B$, $m$-square submatrix $C$ having all enumerating matrix forms in Theorem \ref{the2} have the maximum values on permanent of the matices. Checking $m$-square submatrix $C$, we know that the number of $1$'s in $m$-square submatrix $C$ is $\sigma-(r+1)$.  Now, we consider three subsubcases.

\textbf{Subsubcase 1.2.1}. Assume $\sigma-(r+1)-2m=\sigma-2(m+1)+1-r\equiv0(mod~3)$. Since $\sigma-2(m+1)\equiv1(mod~3)$ and $\sigma-2(m+1)+1-r\equiv0(mod~3)$, we have $2-r=0,-3,-6,\cdots$. This means that $r=2,5,8,\cdots$. By Laplace Expansion Theorem, we get that ${\rm per}B={\rm per}C$. By induction hypothesis, ${\rm per}C=6^{\lfloor\frac{\sigma-2(m+1)-1+2-r}{3} \rfloor} \times2^{\lfloor\frac{3(m+1)-\sigma+r-2}{2} \rfloor}$, and
\begin{eqnarray}\label{equ10}
k_{1}\times6^{\lfloor\frac{\sigma-2(m+1)}{3} \rfloor} \times2^{\lfloor\frac{3(m+1)-\sigma}{2} \rfloor}\leq{\rm per}A={\rm per}B={\rm per}C=6^{\lfloor\frac{\sigma-2(m+1)-1+2-r}{3} \rfloor} \times2^{\lfloor\frac{3(m+1)-\sigma+r-2}{2} \rfloor}.
\end{eqnarray}

\begin{eqnarray}\label{equ11}
{\rm per}C&=&6^{\lfloor\frac{\sigma-2(m+1)-1}{3} \rfloor} \times2^{\lfloor\frac{3(m+1)-\sigma-1}{2} \rfloor}\times6^{\lfloor\frac{2-r}{3} \rfloor} \times2^{\lfloor\frac{r-1}{2} \rfloor}\nonumber\\
&\leq& 6^{\lfloor\frac{\sigma-2(m+1)-1}{3} \rfloor} \times2^{\lfloor\frac{3(m+1)-\sigma}{2} \rfloor}\times6^{\frac{1}{3}}\times\frac{2^{\frac{r-1}{2}}}{6^{\frac{r-1}{3}}}\\
&=&6^{\lfloor\frac{\sigma-2(m+1)-1}{3} \rfloor} \times2^{\lfloor\frac{3(m+1)-\sigma}{2} \rfloor}\times6^{\frac{1}{3}}\times(\frac{\sqrt{2}}{\sqrt[3]{6}})^{r-1}.\nonumber
\end{eqnarray}

Set $r=2$. By (\ref{equ10}), we have ${\rm per}C=6^{\lfloor\frac{\sigma-2(m+1)-1}{3} \rfloor}\times2^{\lfloor\frac{3(m+1)-\sigma}{2} \rfloor}
\leq6^{\lfloor\frac{\sigma-2(m+1)}{3} \rfloor}\times2^{\lfloor\frac{3(m+1)-\sigma}{2} \rfloor}
<k_{1}\cdot6^{\lfloor\frac{\sigma-2(m+1)}{3} \rfloor}\times2^{\lfloor\frac{3(m+1)-\sigma}{2} \rfloor},$
This contradicts to ${\rm per}C\geq k_{1}\times6^{\lfloor\frac{\sigma-2(m+1)}{3} \rfloor} \times2^{\lfloor\frac{3(m+1)-\sigma}{2} \rfloor}$ in (\ref{equ10}). Set $r>2$. By (\ref{equ11}),  we know that if $r$ is enlarged then ${\rm per}C$ is lessened. This indicates ${\rm per}C< k_{1}\times6^{\lfloor\frac{\sigma-2(m+1)}{3} \rfloor} \times2^{\lfloor\frac{3(m+1)-\sigma}{2} \rfloor}$,  a contradiction.

\textbf{Subsubcase 1.2.2}. Assume $\sigma-2(m+1)+1-r\equiv1(mod~3)$. Since $\sigma-2(m+1)\equiv1(mod~3)$ and $\sigma-2(m+1)+1-r\equiv1(mod~3)$, we get that $2-r=1,-2,-5,\cdots$. This means that $r=1,4,7,\cdots$. By Laplace Expansion Theorem, we get that ${\rm per}B={\rm per}C$. By induction hypothesis, there exists $k_{2}\in\{\frac{3}{2},2\}\subseteq \{1, \frac{3}{2}, 2, 3,  4\}$ such that ${\rm per}C= k_{2}\times6^{\lfloor\frac{\sigma-2(m+1)-1+2-r}{3} \rfloor} \times2^{\lfloor\frac{3(m+1)-\sigma+r-2}{2} \rfloor}$, and
\begin{eqnarray}\label{equ12}
k_{1}\times6^{\lfloor\frac{\sigma-2(m+1)}{3} \rfloor} \times2^{\lfloor\frac{3(m+1)-\sigma}{2} \rfloor}\leq{\rm per}A={\rm per}B={\rm per}C\nonumber\\
=k_{2}\times6^{\lfloor\frac{\sigma-2(m+1)-1+2-r}{3} \rfloor} \times2^{\lfloor\frac{3(m+1)-\sigma+r-2}{2} \rfloor}.
\end{eqnarray}
Set $r=1$. By (\ref{equ12}), we have ${\rm per}C=k_{2}\times6^{\lfloor\frac{\sigma-2(m+1)}{3} \rfloor} \times2^{\lfloor\frac{3(m+1)-\sigma-2}{2} \rfloor}
=k_{2}\times6^{\lfloor\frac{\sigma-2(m+1)}{3} \rfloor} \times2^{\lfloor\frac{3(m+1)-\sigma}{2} \rfloor\times\frac{1}{2}}
\leq6^{\lfloor\frac{\sigma-2(m+1)}{3} \rfloor} \times2^{\lfloor\frac{3(m+1)-\sigma}{2} \rfloor}
<k_{1}\times6^{\lfloor\frac{\sigma-2(m+1)}{3} \rfloor} \times2^{\lfloor\frac{3(m+1)-\sigma}{2} \rfloor}.$
This contradicts to ${\rm per}C\geq k_{1}\times6^{\lfloor\frac{\sigma-2(m+1)}{3} \rfloor} \times2^{\lfloor\frac{3(m+1)-\sigma}{2} \rfloor}$ in (\ref{equ12}).
Set $r=4$, By (\ref{equ12}), we have ${\rm per}A={\rm per}B={\rm per}C\leq k_{2}\times6^{\lfloor\frac{\sigma-2(m+1)-3}{3} \rfloor} \times2^{\lfloor\frac{3(m+1)-\sigma+2}{2} \rfloor}\leq\frac{k_{2}}{3}\times6^{\lfloor\frac{\sigma-2(m+1)}{3} \rfloor} \times2^{\lfloor\frac{3(m+1)-\sigma}{2} \rfloor}\leq\frac{2}{3}\times6^{\lfloor\frac{\sigma-2(m+1)}{3} \rfloor} \times2^{\lfloor\frac{3(m+1)-\sigma}{2} \rfloor}<\frac{3}{2}\times6^{\lfloor\frac{\sigma-2(m+1)}{3} \rfloor} \times2^{\lfloor\frac{3(m+1)-\sigma}{2} \rfloor},$
This contradicts to ${\rm per}C\geq k_{1}\times6^{\lfloor\frac{\sigma-2(m+1)}{3} \rfloor} \times2^{\lfloor\frac{3(m+1)-\sigma}{2} \rfloor}$ in (\ref{equ12}). Set $r>4$. By (\ref{equ11}),  we know that if $r$ is enlarged then ${\rm per}C$ is lessened. This indicates ${\rm per}C< k_{1}\times6^{\lfloor\frac{\sigma-2(m+1)}{3} \rfloor} \times2^{\lfloor\frac{3(m+1)-\sigma}{2} \rfloor}$, a contradiction.

\textbf{Subsubcase 1.2.3}. Assume $\sigma-2(m+1)+1-r\equiv2(mod~3)$. Since $\sigma-2(m+1)\equiv1(mod~3)$ and $\sigma-2(m+1)+1-r\equiv2(mod~3)$, we get that $2-r=2,-1,-4,\cdots$. This implies that $r=0,3,6,\cdots$. By Laplace Expansion Theorem, we get that ${\rm per}B={\rm per}C$. By induction hypothesis, there exists $k_{2}\in\{3,4\}\subseteq \{1, \frac{3}{2}, 2, 3,  4\}$ such that ${\rm per}C= k_{2}\times6^{\lfloor\frac{\sigma-2(m+1)-1+2-r}{3} \rfloor} \times2^{\lfloor\frac{3(m+1)-\sigma+r-2}{2} \rfloor}$, and
\begin{eqnarray}\label{equ13}
k_{1}\times6^{\lfloor\frac{\sigma-2(m+1)}{3} \rfloor} \times2^{\lfloor\frac{3(m+1)-\sigma}{2} \rfloor}\leq{\rm per}A={\rm per}B={\rm per}C\nonumber\\
\leq k_{2}\times6^{\lfloor\frac{\sigma-2(m+1)-1+2-r}{3} \rfloor} \times2^{\lfloor\frac{3(m+1)-\sigma+r-2}{2} \rfloor}.
\end{eqnarray}
Set $r=0$. By (\ref{equ12}), we have \begin{eqnarray}\label{equ13x}
k_{1}\times6^{\lfloor\frac{\sigma-2(m+1)}{3} \rfloor} \times2^{\lfloor\frac{3(m+1)-\sigma}{2} \rfloor}
&\leq&{\rm per}A={\rm per}B={\rm per}C \nonumber\\
&=&k_{2}\times6^{\lfloor\frac{\sigma-2(m+1)+1}{3} \rfloor} \times2^{\lfloor\frac{3(m+1)-\sigma-2}{2} \rfloor}\nonumber\\
&=&k_{2}\times6^{\lfloor\frac{\sigma-2(m+1)}{3} \rfloor} \times2^{\lfloor\frac{3(m+1)-\sigma}{2} \rfloor}\times\frac{1}{2}.
\end{eqnarray}
By (\ref{equ13x}), $k_{1}\in\{\frac{3}{2}, 2\}$ and $k_{2}\in\{3, 4\}$, we get that $k_{1}=\frac{3}{2}$ and $k_{2}=3$, or $k_{1}=2$ and $k_{2}=4$.
 If $k_{1}=\frac{3}{2}$ and $k_{2}=3$, then $m$-square submatrix $C$ satisfies item (III) in Theorem \ref{the2}. This implies that $m$-square submatrix $C$ is combinatorially equivalent to $M_{m,\sigma}$. Furthermore,  $(m+1)$-square matrix $A$ is combinatorially equivalent to $S_{(m+1),\sigma}$. If $k_{1}=2$ and $k_{2}=4$, then $m$-square submatrix  $C$ satisfies item (III) in Theorem \ref{the2}. This implies that $m$-square submatrix $C$ is combinatorially equivalent to $N_{m,\sigma}$. Furthermore, $(m+1)$-square matrix $A$ is combinatorially equivalent to $T_{(m+1),\sigma}$. Set $r>0$.  By (\ref{equ11}),  we know that if $r$ is enlarged then ${\rm per}C$ is lessened. This indicates ${\rm per}C< k_{1}\times6^{\lfloor\frac{\sigma-2(m+1)}{3} \rfloor} \times2^{\lfloor\frac{3(m+1)-\sigma}{2} \rfloor}$, a contradiction.

\textbf{Subcase 1.3}. Let $\sigma-2(m+1)\equiv2(mod~3)$.  By simple calculation, there exists $k_{1}\in\{3,4\}$ such  that the permanents of $\{ M_{(m+1),\sigma}, N_{(m+1),\sigma}\}$ equal to $k_{1}\times6^{\lfloor\frac{\sigma-2(m+1)}{3} \rfloor} \times2^{\lfloor\frac{3(m+1)-\sigma}{2} \rfloor}$. This implies that $\mu((m+1),\tau)\geq k_{1}\times6^{\lfloor\frac{\sigma-2(m+1)}{3} \rfloor} \times2^{\lfloor\frac{3(m+1)-\sigma}{2} \rfloor}$. By induction hypothesis, for any $m$-square submatrix $C$ of $B$, $m$-square submatrix $C$ having all enumerating matrix forms in Theorem \ref{the2} has the maximum values on permanent of the matrices. Checking $m$-square submatrix $C$, we know that the number of $1$'s in $m$-square submatrix $C$ is $\sigma-(r+1)$.  Now, we consider three subsubcases.

\textbf{Subsubcase 1.3.1}. Assume $\sigma-(r+1)-2((m+1)-1)=\sigma-2(m+1)+1-r\equiv0(mod~3)$. Since $\sigma-2(m+1)\equiv2(mod~3)$ and $\sigma-2(m+1)+1-r\equiv0(mod~3)$, we get that $3-r=0,-3,-6,\cdots$. This means that $r=3,6,9,\cdots$. By Laplace Expansion Theorem, we get that ${\rm per}B={\rm per}C$. By induction hypothesis, ${\rm per}C= 6^{\lfloor\frac{\sigma-2(m+1)-1+2-r}{3} \rfloor} \times2^{\lfloor\frac{3(m+1)-\sigma+r-2}{2} \rfloor}$. So,
\begin{eqnarray}\label{equ14}
k_{1}\times6^{\lfloor\frac{\sigma-2(m+1)}{3} \rfloor} \times2^{\lfloor\frac{3(m+1)-\sigma}{2} \rfloor}\leq{\rm per}A={\rm per}B={\rm per}C=6^{\lfloor\frac{\sigma-2(m+1)-1+2-r}{3} \rfloor} \times2^{\lfloor\frac{3(m+1)-\sigma+r-2}{2} \rfloor}.
\end{eqnarray}
By (\ref{equ14}), we have
\begin{eqnarray}\label{equ15}
{\rm per}C&=&k_{1}\times6^{\lfloor\frac{\sigma-2(m+1)-2}{3} \rfloor} \times2^{\lfloor\frac{3(m+1)-\sigma-1}{2} \rfloor}\times6^{\lfloor\frac{3-r}{3} \rfloor} \times2^{\lfloor\frac{r-1}{2} \rfloor}\nonumber\\
&\leq& k_{1}\times6^{\lfloor\frac{\sigma-2(m+1)-2}{3} \rfloor} \times2^{\lfloor\frac{3(m+1)-\sigma}{2} \rfloor}\times6^{\frac{2}{3}}\times\frac{2^{\frac{r-1}{2}}}{6^{\frac{r-1}{3}}}\\
&=&k_{1}\times6^{\lfloor\frac{\sigma-2(m+1)-2}{3} \rfloor} \times2^{\lfloor\frac{3(m+1)-\sigma}{2} \rfloor}\times6^{\frac{2}{3}}\times(\frac{\sqrt{2}}{\sqrt[3]{6}})^{r-1}.\nonumber
\end{eqnarray}
Set $r=3$. By (\ref{equ14}), we have ${\rm per}C=6^{\lfloor\frac{\sigma-2(m+1)-2}{3} \rfloor}\times2^{\lfloor\frac{3(m+1)-\sigma+1}{2} \rfloor}
\leq2\times6^{\lfloor\frac{\sigma-2(m+1)}{3} \rfloor}\times2^{\lfloor\frac{3(m+1)-\sigma}{2} \rfloor}
<k_{1}\times6^{\lfloor\frac{\sigma-2(m+1)}{3} \rfloor}\times2^{\lfloor\frac{3(m+1)-\sigma}{2} \rfloor}.$
 This contradicts to ${\rm per}C\geq k_{1}\times6^{\lfloor\frac{\sigma-2(m+1)}{3} \rfloor} \times2^{\lfloor\frac{3(m+1)-\sigma}{2} \rfloor}$ in (\ref{equ14}). Set  $r>3$.  By (\ref{equ11}),  we know that if $r$ is enlarged then ${\rm per}C$ is lessened. This indicates ${\rm per}C< k_{1}\times6^{\lfloor\frac{\sigma-2(m+1)}{3} \rfloor} \times2^{\lfloor\frac{3(m+1)-\sigma}{2} \rfloor}$, a contradiction.

\textbf{Subsubcase 1.3.2}. Assume $\sigma-2(m+1)+1-r\equiv1(mod~3)$. Since $\sigma-2(m+1)\equiv2(mod~3)$ and $\sigma-2(m+1)+1-r\equiv1(mod~3)$, we get that $3-r=1,-2,-5,\cdots$. This means that $r=2,5,8,\cdots$. By Laplace Expansion Theorem, we get that ${\rm per}B={\rm per}C$. By induction hypothesis, there exists $k_{2}\in\{\frac{3}{2},2\}$ such that ${\rm per}C=k_{2}\times6^{\lfloor\frac{\sigma-2(m+1)-1+2-r}{3} \rfloor} \times2^{\lfloor\frac{3(m+1)-\sigma+r-2}{2} \rfloor}$, and
\begin{eqnarray}\label{equ16}
k_{1}\times6^{\lfloor\frac{\sigma-2(m+1)}{3} \rfloor} \times2^{\lfloor\frac{3(m+1)-\sigma}{2} \rfloor}\leq{\rm per}A={\rm per}B={\rm per}C=\nonumber\\
k_{2}\times6^{\lfloor\frac{\sigma-2(m+1)-1+2-r}{3} \rfloor} \times2^{\lfloor\frac{3(m+1)-\sigma+r-2}{2} \rfloor}.
\end{eqnarray}
 Set $r=2$. If $k_{2}=\frac{3}{2}$, then  we know that $m$-square submatrix $C$ satisfies item (II) in Theorem \ref{the2}. This implies that $3(m+1)-\sigma-1$ is even and $3(m+1)-\sigma$ is odd.  By (\ref{equ16}), we have ${\rm per}C=\frac{3}{2}\times6^{\lfloor\frac{\sigma-2(m+1)-1}{3} \rfloor} \times2^{\lfloor\frac{3(m+1)-\sigma-1}{2} \rfloor}
=\frac{3}{2}\times6^{\lfloor\frac{\sigma-2(m+1)}{3} \rfloor} \times2^{\lfloor\frac{3(m+1)-\sigma}{2} \rfloor}
<k_{1}\times6^{\lfloor\frac{\sigma-2(m+1)}{3} \rfloor} \times2^{\lfloor\frac{3(m+1)-\sigma}{2} \rfloor}.$
 This contradicts to ${\rm per}C\geq k_{1}\times6^{\lfloor\frac{\sigma-2(m+1)}{3} \rfloor} \times2^{\lfloor\frac{3(m+1)-\sigma}{2} \rfloor}$ in (\ref{equ16}).  If $k_{2}=2$, then  we know that $m$-square submatrix $C$ satisfies item (II) in Theorem \ref{the2}. This implies that $3(m+1)-\sigma-1$ is odd and $3(m+1)-\sigma$ is odd.  By (\ref{equ16}), we have ${\rm per}C=2\times6^{\lfloor\frac{\sigma-2(m+1)-1}{3}\rfloor} \times2^{\lfloor\frac{3(m+1)-\sigma-1}{2} \rfloor}=2\times6^{\lfloor\frac{\sigma-2(m+1)}{3} \rfloor} \times2^{\lfloor\frac{3(m+1)-\sigma}{2} \rfloor}\times\frac{1}{2}=6^{\lfloor\frac{\sigma-2(m+1)}{3} \rfloor} \times2^{\lfloor\frac{3(m+1)-\sigma}{2} \rfloor}<k_{1}\times6^{\lfloor\frac{\sigma-2(m+1)}{3} \rfloor} \times2^{\lfloor\frac{3(m+1)-\sigma}{2} \rfloor}.$
This contradicts to ${\rm per}C\geq k_{1}\times6^{\lfloor\frac{\sigma-2(m+1)}{3} \rfloor} \times2^{\lfloor\frac{3(m+1)-\sigma}{2} \rfloor}$ in (\ref{equ16}). Set $r>2$.  By (\ref{equ11}),  we know that if $r$ is enlarged then ${\rm per}C$ is lessened. This indicates ${\rm per}C< k_{1}\times6^{\lfloor\frac{\sigma-2(m+1)}{3} \rfloor} \times2^{\lfloor\frac{3(m+1)-\sigma}{2} \rfloor}$, a contradiction.

\textbf{Subsubcase 1.3.3}. Assume $\sigma-2(m+1)+1-r\equiv2(mod~3)$.  Since $\sigma-2(m+1)\equiv2(mod~3)$ and $\sigma-2(m+1)+1-r\equiv2(mod~3)$, we have $3-r=2,-1,-4,\cdots$. This means that $r=1,4,7,\cdots$. By Laplace Expansion Theorem, we get that ${\rm per}B={\rm per}C$. By induction hypothesis, there exists $k_{2}\in\{3,4\}$ such that ${\rm per}C= k_{2}\times6^{\lfloor\frac{\sigma-2(m+1)-1+2-r}{3} \rfloor} \times2^{\lfloor\frac{3(m+1)-\sigma+r-2}{2} \rfloor}$, and
\begin{eqnarray}\label{equ17}
k_{1}\times6^{\lfloor\frac{\sigma-2(m+1)}{3} \rfloor} \times2^{\lfloor\frac{3(m+1)-\sigma}{2} \rfloor}\leq{\rm per}A={\rm per}B={\rm per}C=\nonumber\\
k_{2}\times6^{\lfloor\frac{\sigma-2(m+1)-1+2-r}{3} \rfloor} \times2^{\lfloor\frac{3(m+1)-\sigma+r-2}{2} \rfloor}.
\end{eqnarray}
Set $r=1$. If  $k_{2}=4$, then  we know that $m$-square submatrix $C$ satisfies item (III) in Theorem \ref{the2}. This implies that  $3(m+1)-\sigma-1$ is odd  and $3(m+1)-\sigma$ is even. By (\ref{equ17}), we have
${\rm per}C=4\times6^{\lfloor\frac{\sigma-2(m+1)}{3} \rfloor} \times2^{\lfloor\frac{3(m+1)-\sigma-2}{2} \rfloor}\times\frac{1}{2}
=2\times6^{\lfloor\frac{\sigma-2(m+1)}{3} \rfloor}\times2^{\lfloor\frac{3(m+1)-\sigma}{2} \rfloor}
<k_{1}\times6^{\lfloor\frac{\sigma-2(m+1)}{3} \rfloor}\times2^{\lfloor\frac{3(m+1)-\sigma}{2} \rfloor},$
This contradicts to ${\rm per}C\geq k_{1}\times6^{\lfloor\frac{\sigma-2(m+1)}{3} \rfloor} \times2^{\lfloor\frac{3(m+1)-\sigma}{2} \rfloor}$ in (\ref{equ17}). If $k_{2}=3$, then  we know that $m$-square  submatrix $C$ satisfies item (III) in Theorem \ref{the2}. It is easy to obtain that ${\rm per}C<k_{1}\times6^{\lfloor\frac{\sigma-2(m+1)}{3} \rfloor}\times2^{\lfloor\frac{3(m+1)-\sigma}{2} \rfloor}.$  This contradicts to ${\rm per}C\geq k_{1}\times6^{\lfloor\frac{\sigma-2(m+1)}{3} \rfloor} \times2^{\lfloor\frac{3(m+1)-\sigma}{2} \rfloor}$ in (\ref{equ17}).
Set $r>1$. By (\ref{equ11}),  we know that if $r$ is enlarged then ${\rm per}C$ is lessened. This indicates ${\rm per}C< k_{1}\times6^{\lfloor\frac{\sigma-2(m+1)}{3} \rfloor} \times2^{\lfloor\frac{3(m+1)-\sigma}{2} \rfloor}$, a contradiction.


Now we determine  the form of $(m+1)$-square matrix $D$. Assume that $m$-square submatrices $D_{1}$ and $D_{2}$ have $s$~$1$'s and $t$~$1$'s, respectively. By induction hypothesis and simple calculating,  ${\rm per}D(1,1)=k_{2}\times6^{\lfloor\frac{\sigma-2(m+1)-s}{3} \rfloor} \times2^{\lfloor\frac{3(m+1)-\sigma+s-1}{2} \rfloor}$ and   ${\rm per}D(1,2)=k_{3}\times6^{\lfloor\frac{\sigma-2(m+1)-t}{3} \rfloor} \times2^{\lfloor\frac{3(m+1)-\sigma+t-1}{2} \rfloor}$, where $k_{2}, k_{3}\in \{1, \frac{3}{2}, 2, 3,  4\} $. In particular, if $n=m+1$ then the permanents of
$F_{m+1,\sigma}$, $F_{m+1,\sigma}^{*}$, $U_{m+1,\sigma}$, $V_{m+1,\sigma}$, $W_{m+1,\sigma}$, $X_{m+1,\sigma}$, $Y_{m+1,\sigma}$, $Y_{m+1,\sigma}$, $F_{m+1,\sigma}^{*}$, $R_{m+1,\sigma}$, $S_{m+1,\sigma}$, $T_{m+1,\sigma}$,  $Q_{m+1,\sigma}$, $P_{m+1,\sigma}^{*}$, $M_{m+1,\sigma}$ and $N_{m+1,\sigma}$ are equal to  $k_{1}\times6^{\lfloor\frac{\sigma-2(m+1)}{3} \rfloor} \times2^{\lfloor\frac{3(m+1)-\sigma}{2} \rfloor}$,  where $k_{1}\in\{1, \frac{3}{2}, 2, 3,  4\}$.
These imply that ${\rm per}A=\mu(m+1,\tau)\geq k_{1}\times6^{\lfloor\frac{\sigma-2(m+1)}{3} \rfloor} \times2^{\lfloor\frac{3(m+1)-\sigma}{2} \rfloor}$, where $k_{1}\in\{1, \frac{3}{2}, 2, 3,  4\}$. By Laplace Expansion Theorem, we get that ${\rm per}D={\rm per}D(1,1)+{\rm per}D(1,2)$. Checking $m$-square submatrix $D(1,1)$, we know that the number of $1$'s in $m$-square submatrix $D(1,1)$ is $\sigma-s-2$.
 By induction hypothesis, ${\rm per}D= k_{2}\times6^{\lfloor\frac{\sigma-2(m+1)-s}{3} \rfloor} \times2^{\lfloor\frac{3(m+1)-\sigma+s-1}{2} \rfloor}+k_{3}\times6^{\lfloor\frac{\sigma-2(m+1)-t}{3} \rfloor} \times2^{\lfloor\frac{3(m+1)-\sigma+t-1}{2} \rfloor}$, and
\begin{eqnarray}\label{equ18}
k_{1}\times6^{\lfloor\frac{\sigma-2(m+1)}{3} \rfloor} \times2^{\lfloor\frac{3(m+1)-\sigma}{2} \rfloor}\leq{\rm per}A={\rm per}D=\nonumber\\
k_{2}\times6^{\lfloor\frac{\sigma-s-2-2m}{3} \rfloor} \times2^{\lfloor\frac{3m-(\sigma-s-2)}{2} \rfloor}+k_{3}\times6^{\lfloor\frac{\sigma-t-2-2m}{3} \rfloor} \times2^{\lfloor\frac{3m-(\sigma-t-2)}{2} \rfloor}\\
=k_{2}\times6^{\lfloor\frac{\sigma-2(m+1)-s}{3} \rfloor} \times2^{\lfloor\frac{3(m+1)-\sigma+s-1}{2}
\rfloor}+k_{3}\times6^{\lfloor\frac{\sigma-2(m+1)-t}{3} \rfloor} \times2^{\lfloor\frac{3(m+1)-\sigma+t-1}{2} \rfloor}.\nonumber
\end{eqnarray}

For the purpose of further proof, we give a claim as follows.

\begin{claim}\label{cla2}
(i) When $\sigma-2(m+1)\equiv2(mod~3)$ and $3(m+1)-\sigma$ is odd.  If $s$ is enlarged,  then $k_{2}\times6^{\lfloor\frac{\sigma-2(m+1)-s}{3} \rfloor} \times2^{\lfloor\frac{3(m+1)-\sigma+s-1}{2} \rfloor}$ is lessened.

(ii) When $\sigma-2(m+1)\equiv2(mod~3)$ and $3(m+1)-\sigma$ is even. If $s$ is enlarged,  then  $k_{2}\times6^{\lfloor\frac{\sigma-2(m+1)-s}{3} \rfloor} \times2^{\lfloor\frac{3(m+1)-\sigma+s-1}{2} \rfloor}$ is lessened.

(iii) When $\sigma-2(m+1)\equiv1(mod~3)$ and $3(m+1)-\sigma$ is odd. If $s$ is enlarged,  then $k_{2}\times6^{\lfloor\frac{\sigma-2(m+1)-s}{3} \rfloor} \times2^{\lfloor\frac{3(m+1)-\sigma+s-1}{2} \rfloor}$ is lessened.

(iv) When $\sigma-2(m+1)\equiv1(mod~3)$ and $3(m+1)-\sigma$ is even. If $s$ is enlarged,  then $k_{2}\times6^{\lfloor\frac{\sigma-2(m+1)-s}{3} \rfloor} \times2^{\lfloor\frac{3(m+1)-\sigma+s-1}{2} \rfloor}$ is lessened.

(v) When $\sigma-2(m+1)\equiv0(mod~3)$.  If $s$ is enlarged,  then $k_{2}\times6^{\lfloor\frac{\sigma-2(m+1)-s}{3} \rfloor} \times2^{\lfloor\frac{3(m+1)-\sigma+s-1}{2} \rfloor}$ is lessened.
\end{claim}

\textbf{Proof of Claim \ref{cla2}.} (i) Since  $\sigma-2(m+1)\equiv2(mod~3)$ and $3(m+1)-\sigma$ is odd, we get that $k_{2}\times6^{\lfloor\frac{\sigma-2(m+1)-s}{3} \rfloor} \times2^{\lfloor\frac{3(m+1)-\sigma+s-1}{2} \rfloor}=k_{2}\times6^{\lfloor\frac{\sigma-2(m+1)-s}{3} \rfloor} \times2^{\lfloor\frac{3(m+1)-\sigma+s-1}{2} \rfloor}=k_{2}\times6^{\lfloor\frac{\sigma-2(m+1)-2}{3} \rfloor} \times2^{\lfloor\frac{3(m+1)-\sigma-1}{2} \rfloor}\times6^{\lfloor\frac{2-s}{3} \rfloor} \times2^{\lfloor\frac{s}{2} \rfloor}<k_{2}\times6^{\lfloor\frac{\sigma-2(m+1)-2}{3} \rfloor} \times2^{\lfloor\frac{3(m+1)-\sigma-1}{2} \rfloor}\times6^{\frac{3-s}{3}} \times2^{\frac{s}{2}}=k_{2}\times6\times6^{\lfloor\frac{\sigma-2(m+1)-2}{3} \rfloor} \times2^{\lfloor\frac{3(m+1)-\sigma-1}{2} \rfloor}\times(\frac{\sqrt{2}}{\sqrt[3]{6}})^{s}$.  Since $\frac{\sqrt{2}}{\sqrt[3]{6}}< 1$,
we know that if $s$ is enlarged then $k_{2}\times6^{\lfloor\frac{\sigma-2(m+1)-s}{3} \rfloor} \times2^{\lfloor\frac{3(m+1)-\sigma+s-1}{2} \rfloor}$ is lessened.

(ii) Since $\sigma-2(m+1)\equiv2(mod~3)$ and $3(m+1)-\sigma$ is even, we obtain that $k_{2}\times6^{\lfloor\frac{\sigma-2(m+1)-s}{3} \rfloor} \times2^{\lfloor\frac{3(m+1)-\sigma+s-1}{2} \rfloor}=k_{2}\times6^{\lfloor\frac{\sigma-2(m+1)-2-s+2}{3} \rfloor} \times2^{\lfloor\frac{3(m+1)-\sigma+s-1}{2} \rfloor}=k_{2}\times6^{\lfloor\frac{\sigma-2(m+1)-2}{3} \rfloor} \times2^{\lfloor\frac{3(m+1)-\sigma}{2} \rfloor}\times6^{\lfloor\frac{2-s}{3} \rfloor} \times2^{\lfloor\frac{s-1}{2} \rfloor}\leq k_{2}\times6\times6^{\lfloor\frac{\sigma-2(m+1)-2}{3} \rfloor} \times2^{\lfloor\frac{3(m+1)-\sigma}{2} \rfloor}\times6^{\lfloor\frac{-s}{3} \rfloor} \times2^{\lfloor\frac{s}{2} \rfloor}<k_{2}\times6\times6^{\lfloor\frac{\sigma-2(m+1)-2}{3} \rfloor} \times2^{\lfloor\frac{3(m+1)-\sigma-1}{2} \rfloor}\times(\frac{\sqrt{2}}{\sqrt[3]{6}})^{s}$. Since $\frac{\sqrt{2}}{\sqrt[3]{6}}< 1$, we know that if $s$ is enlarged then $k_{2}\times6^{\lfloor\frac{\sigma-2(m+1)-s}{3} \rfloor} \times2^{\lfloor\frac{3(m+1)-\sigma+s-1}{2} \rfloor}$ is lessened.

(iii) Since $\sigma-2(m+1)\equiv1(mod~3)$ and $3(m+1)-\sigma$ is even, we get that $k_{2}\times6^{\lfloor\frac{\sigma-2(m+1)-s}{3} \rfloor} \times2^{\lfloor\frac{3(m+1)-\sigma+s-1}{2} \rfloor}=k_{2}\times6^{\lfloor\frac{\sigma-2(m+1)-2-s+2}{3} \rfloor} \times2^{\lfloor\frac{3(m+1)-\sigma+s-1}{2} \rfloor}=k_{2}\times6^{\lfloor\frac{\sigma-2(m+1)-2}{3} \rfloor} \times2^{\lfloor\frac{3(m+1)-\sigma}{2} \rfloor}\times6^{\lfloor\frac{2-s}{3} \rfloor} \times2^{\lfloor\frac{s-1}{2} \rfloor}\leq k_{2}\times6\times6^{\lfloor\frac{\sigma-2(m+1)-2}{3} \rfloor} \times2^{\lfloor\frac{3(m+1)-\sigma}{2} \rfloor}\times6^{\lfloor\frac{-s}{3} \rfloor} \times2^{\lfloor\frac{s}{2} \rfloor}<k_{2}\times6\times6^{\lfloor\frac{\sigma-2(m+1)-2}{3} \rfloor} \times2^{\lfloor\frac{3(m+1)-\sigma-1}{2} \rfloor}\times(\frac{\sqrt{2}}{\sqrt[3]{6}})^{s}$. Since $\frac{\sqrt{2}}{\sqrt[3]{6}}< 1$, we know that if $s$ is enlarged then $k_{2}\times6^{\lfloor\frac{\sigma-2(m+1)-s}{3} \rfloor} \times2^{\lfloor\frac{3(m+1)-\sigma+s-1}{2} \rfloor}$ is lessened.

(iv) Since $\sigma-2(m+1)\equiv1(mod~3)$ and $3(m+1)-\sigma$ is even, we get that $k_{2}\times6^{\lfloor\frac{\sigma-2(m+1)-s}{3} \rfloor} \times2^{\lfloor\frac{3(m+1)-\sigma+s-1}{2} \rfloor}=k_{2}\times6^{\lfloor\frac{\sigma-2(m+1)-1-s+1}{3} \rfloor} \times2^{\lfloor\frac{3(m+1)-\sigma+s-1}{2} \rfloor}=k_{2}\times6^{\lfloor\frac{\sigma-2(m+1)-1}{3} \rfloor} \times2^{\lfloor\frac{3(m+1)-\sigma}{2} \rfloor}\times6^{\lfloor\frac{1-s}{3} \rfloor} \times2^{\lfloor\frac{s-1}{2} \rfloor}\leq k_{2}\times6\times6^{\lfloor\frac{\sigma-2(m+1)-1}{3} \rfloor} \times2^{\lfloor\frac{3(m+1)-\sigma}{2} \rfloor}\times6^{\lfloor\frac{-s}{3} \rfloor} \times2^{\lfloor\frac{s}{2} \rfloor}\leq k_{2}\times6\times6^{\lfloor\frac{\sigma-2(m+1)-2}{3} \rfloor} \times2^{\lfloor\frac{3(m+1)-\sigma-1}{2} \rfloor}\times6^{\frac{-s}{3}} \times2^{\frac{s}{2}}=k_{2}\times6\times6^{\lfloor\frac{\sigma-2(m+1)-2}{3} \rfloor} \times2^{\lfloor\frac{3(m+1)-\sigma-1}{2} \rfloor}\times(\frac{\sqrt{2}}{\sqrt[3]{6}})^{s}$. Since $\frac{\sqrt{2}}{\sqrt[3]{6}}< 1$, we know that if $s$ is enlarged then $k_{2}\times6^{\lfloor\frac{\sigma-2(m+1)-s}{3} \rfloor} \times2^{\lfloor\frac{3(m+1)-\sigma+s-1}{2} \rfloor}$ is lessened.

(v) Since $\sigma-2(m+1)\equiv0(mod~3)$, we get that $k_{2}\times6^{\lfloor\frac{\sigma-2(m+1)-s}{3} \rfloor} \times2^{\lfloor\frac{3(m+1)-\sigma+s-1}{2} \rfloor}=k_{2}\times6^{\lfloor\frac{\sigma-2(m+1)-s}{3} \rfloor} \times2^{\lfloor\frac{3(m+1)-\sigma+s-1}{2} \rfloor}=k_{2}\times6^{\lfloor\frac{\sigma-2(m+1)}{3} \rfloor} \times2^{\lfloor\frac{3(m+1)-\sigma}{2} \rfloor}\times6^{\lfloor\frac{-s}{3} \rfloor} \times2^{\lfloor\frac{s}{2} \rfloor}\leq k_{2}\times6^{\lfloor\frac{\sigma-2(m+1)-1}{3} \rfloor} \times2^{\lfloor\frac{3(m+1)-\sigma}{2} \rfloor}\times6^{\lfloor\frac{-s}{3} \rfloor} \times2^{\lfloor\frac{s}{2} \rfloor}\leq k_{2}\times6^{\lfloor\frac{\sigma-2(m+1)}{3} \rfloor} \times2^{\lfloor\frac{3(m+1)-\sigma}{2} \rfloor}\times6^{\frac{-s}{3}} \times2^{\frac{s}{2}}=k_{2}\times6^{\lfloor\frac{\sigma-2(m+1)}{3} \rfloor} \times2^{\lfloor\frac{3(m+1)-\sigma}{2} \rfloor}\times(\frac{\sqrt{2}}{\sqrt[3]{6}})^{s}$. Since $\frac{\sqrt{2}}{\sqrt[3]{6}}< 1$, we know that if $s$ is enlarged then $k_{2}\times6^{\lfloor\frac{\sigma-2(m+1)-s}{3} \rfloor} \times2^{\lfloor\frac{3(m+1)-\sigma+s-1}{2} \rfloor}$ is lessened. \qquad\qquad\qquad\qquad\qquad\qquad\qquad\qquad\qquad\qquad\qquad\qquad\qquad\qquad$\Box$

By (\ref{equ18}),  we  consider  $k_{2}=k_{3}$ and $k_{2}\neq k_{3}$. In particular, we need to point out that if $s=0$ in $D_{1}$ or $t=0$ in  $D_{2}$, then the  form of $(m+1)$-square matrix $D$ is equal to the  form of $(m+1)$-square matrix $B$. This case is discussed as above. Nextly, we only consider $s\geq1$ and $t\geq1$ in $D_{1}$ and $D_{2}$.

 Assume $k_{2}=k_{3}$. We will discuss 25 cases, i.e.,  $k_{2}=k_{3}=4$, $\sigma-2(m+1)\equiv2(mod~3)$ and $3(m+1)-\sigma$ is odd;  $k_{2}=k_{3}=4$, $\sigma-2(m+1)\equiv2(mod~3)$ and  $3(m+1)-\sigma$ is odd; $k_{2}=k_{3}=4$, $\sigma-2(m+1)\equiv1(mod~3)$ and $3(m+1)-\sigma$ is odd; $k_{2}=k_{3}=4$,  $\sigma-2(m+1)\equiv1(mod~3)$ and $3(m+1)-\sigma$ is odd; $k_{2}=k_{3}=4$ and  $\sigma-2(m+1)\equiv0(mod~3)$;
 $k_{2}=k_{3}=3$,   $\sigma-2(m+1)\equiv2(mod~3)$ and $3(m+1)-\sigma$ is odd;  $k_{2}=k_{3}=3$,  $\sigma-2(m+1)\equiv2(mod~3)$ and $3(m+1)-\sigma$ is odd; $k_{2}=k_{3}=3$,  $\sigma-2(m+1)\equiv1(mod~3)$ and $3(m+1)-\sigma$ is odd; $k_{2}=k_{3}=3$, $\sigma-2(m+1)\equiv1(mod~3)$ and $3(m+1)-\sigma$ is even; $k_{2}=k_{3}=4$ and  $\sigma-2(m+1)\equiv0(mod~3)$;
 $k_{2}=k_{3}=2$,  $\sigma-2(m+1)\equiv2(mod~3)$ and $3(m+1)-\sigma$ is odd;  $k_{2}=k_{3}=2$,  $\sigma-2(m+1)\equiv2(mod~3)$ and $3(m+1)-\sigma$ is even; $k_{2}=k_{3}=2$,  $\sigma-2(m+1)\equiv1(mod~3)$ and $3(m+1)-\sigma$ is odd; $k_{2}=k_{3}=2$,  $\sigma-2(m+1)\equiv1(mod~3)$ and $3(m+1)-\sigma$ is even; $k_{2}=k_{3}=2$ and   $\sigma-2(m+1)\equiv1(mod~3)$;
 $k_{2}=k_{3}=\frac{3}{2}$,  $\sigma-2(m+1)\equiv2(mod~3)$ and $3(m+1)-\sigma$ is odd;  $k_{2}=k_{3}=\frac{3}{2}$,  $\sigma-2(m+1)\equiv2(mod~3)$ and $3(m+1)-\sigma$ is even; $k_{2}=k_{3}=\frac{3}{2}$,  $\sigma-2(m+1)\equiv1(mod~3)$ and $3(m+1)-\sigma$ is odd; $k_{2}=k_{3}=\frac{3}{2}$,  $\sigma-2(m+1)\equiv1(mod~3)$ and $3(m+1)-\sigma$ is even; $k_{2}=k_{3}=\frac{3}{2}$ and  $\sigma-2(m+1)\equiv0(mod~3)$;
 $k_{2}=k_{3}=1$,  $\sigma-2(m+1)\equiv2(mod~3)$ and $3(m+1)-\sigma$ is odd;  $k_{2}=k_{3}=1$,  $\sigma-2(m+1)\equiv2(mod~3)$ and $3(m+1)-\sigma$ is even; $k_{2}=k_{3}=1$, $\sigma-2(m+1)\equiv1(mod~3)$ and $3(m+1)-\sigma$ is odd; $k_{2}=k_{3}=1$, $\sigma-2(m+1)\equiv1(mod~3)$ and $3(m+1)-\sigma$ is even; $k_{2}=k_{3}=1$ and  $\sigma-2(m+1)\equiv0(mod~3)$. The specific discussion is as follows.

Assume that $k_{2}=k_{3}=4$ and $\sigma-2(m+1)\equiv2(mod~3)$, and $3(m+1)-\sigma$ is odd. This implies that $\sigma-2(m+1)-s\equiv2~(mod~3)$, and $3(m+1)-\sigma+s-1$ is odd. Since $\sigma-2(m+1)-s\equiv2~(mod~3)$, and $3(m+1)-\sigma+s-1$ is odd, we have $s=3,6,9,\cdots$, and $s-1=0,2,4,\cdots$. This means that $s=3,9,15,\cdots$. Similarly, we have $t=3,9,15,\cdots$.
Set that $s=3$ and $t=3$.   By (\ref{equ18}), we have ${\rm per}A={\rm per}D=4\times6^{\lfloor\frac{\sigma-2(m+1)-s}{3} \rfloor} \times2^{\lfloor\frac{3(m+1)-\sigma+s-1}{2} \rfloor}+4\times6^{\lfloor\frac{\sigma-2(m+1)-t}{3} \rfloor} \times2^{\lfloor\frac{3(m+1)-\sigma+t-1}{2} \rfloor}=\frac{4}{3}\times6^{\lfloor\frac{\sigma-2(m+1)}{3} \rfloor} \times2^{\lfloor\frac{3(m+1)-\sigma}{2} \rfloor}<4\times6^{\lfloor\frac{\sigma-2(m+1)}{3} \rfloor} \times2^{\lfloor\frac{3(m+1)-\sigma}{2} \rfloor}=4\times6^{\lfloor\frac{\sigma-2(m+1)}{3} \rfloor} \times2^{\lfloor\frac{3(m+1)-\sigma}{2} \rfloor}$. This contradicts to ${\rm per}D\geq 4\times6^{\lfloor\frac{\sigma-2(m+1)}{3} \rfloor} \times2^{\lfloor\frac{3(m+1)-\sigma}{2} \rfloor}$ in (\ref{equ18}).
Set $s>3$ or $t>3$.  By (i) of Claim \ref{cla2}, we have ${\rm per}D<4\times6^{\lfloor\frac{\sigma-2(m+1)}{3} \rfloor} \times2^{\lfloor\frac{3(m+1)-\sigma}{2} \rfloor}$. This contradicts to ${\rm per}D\geq4\times6^{\lfloor\frac{\sigma-2(m+1)}{3} \rfloor} \times2^{\lfloor\frac{3(m+1)-\sigma}{2} \rfloor}$ in (\ref{equ18}).

Assume that $k_{2}=k_{3}=4$ and $\sigma-2(m+1)\equiv2(mod~3)$, and $3(m+1)-\sigma$ is odd.  This implies that  $\sigma-2(m+1)-s\equiv2(mod~3)$, and $3(m+1)-\sigma+s-1$ is odd. Since $\sigma-2(m+1)-s\equiv2(mod~3)$ and $3(m+1)-\sigma+s-1$ is odd,  we have  $s=3,6,9,\cdots$, and $s-1=1,3,5,\cdots$. This means that $s=6,12,18,\cdots$. Similarly, we have $t=6,12,18,\cdots$.  Set $s=6$ and $t=6$. By (\ref{equ18}), we have  ${\rm per}A={\rm per}D=4\times6^{\lfloor\frac{\sigma-2(m+1)-s}{3} \rfloor} \times2^{\lfloor\frac{3(m+1)-\sigma+s-1}{2} \rfloor}+4\times6^{\lfloor\frac{\sigma-2(m+1)-t}{3} \rfloor} \times2^{\lfloor\frac{3(m+1)-\sigma+t-1}{2} \rfloor}=4\times2\times6^{\lfloor\frac{\sigma-2(m+1)-6}{3} \rfloor} \times2^{\lfloor\frac{3(m+1)-\sigma+5}{2} \rfloor}=4\times2\times4\times\frac{1}{36}\times6^{\lfloor\frac{\sigma-2(m+1)-6}{3} \rfloor} \times2^{\lfloor\frac{3(m+1)-\sigma+5}{2} \rfloor}=\frac{8}{9}\times6^{\lfloor\frac{\sigma-2(m+1)-6}{3} \rfloor} \times2^{\lfloor\frac{3(m+1)-\sigma+5}{2} \rfloor}<3\times6^{\lfloor\frac{\sigma-2(m+1)}{3} \rfloor} \times2^{\lfloor\frac{3(m+1)-\sigma}{2} \rfloor}=3\times6^{\lfloor\frac{\sigma-2(m+1)}{3} \rfloor}\times2^{\lfloor\frac{3(m+1)-\sigma}{2} \rfloor}$. This contradicts to ${\rm per}D\geq 3\times6^{\lfloor\frac{\sigma-2(m+1)}{3} \rfloor} \times2^{\lfloor\frac{3(m+1)-\sigma}{2} \rfloor}$ in (\ref{equ18}). Set that $s>6$ or $t>6$.  By (ii) of Claim \ref{cla2}, we have ${\rm per}D<3\times6^{\lfloor\frac{\sigma-2(m+1)}{3} \rfloor} \times2^{\lfloor\frac{3(m+1)-\sigma}{2} \rfloor}$. This contradicts to ${\rm per}D\geq 3\times6^{\lfloor\frac{\sigma-2(m+1)}{3} \rfloor} \times2^{\lfloor\frac{3(m+1)-\sigma}{2} \rfloor}$ in (\ref{equ18}).

Assume that $k_{2}=k_{3}=4$ and  $\sigma-2(m+1)\equiv1(mod~3)$, and $3(m+1)-\sigma$ is odd. This implies that  $\sigma-2(m+1)-s\equiv2(mod~3)$ and $3(m+1)-\sigma+s-1$ is odd. Since $\sigma-2(m+1)-s\equiv2(mod~3)$ and $3(m+1)-\sigma+s-1$ is odd, we have  $s=2,5,8,\cdots$, and $s-1=0,2,4,\cdots$. This means that  $s=5,11,17,\cdots$. Similarly, we have $t=5,11,17,\cdots$. Set that $s=5$ and $t=5$. By (\ref{equ18}), we have ${\rm per}A={\rm per}D=4\times6^{\lfloor\frac{\sigma-2(m+1)-s}{3} \rfloor} \times2^{\lfloor\frac{3(m+1)-\sigma+s-1}{2} \rfloor}+4\times6^{\lfloor\frac{\sigma-2(m+1)-t}{3} \rfloor} \times2^{\lfloor\frac{3(m+1)-\sigma+t-1}{2} \rfloor}=\frac{8}{9}\times6^{\lfloor\frac{\sigma-2(m+1)}{3} \rfloor} \times2^{\lfloor\frac{3(m+1)-\sigma}{2} \rfloor}<2\times6^{\lfloor\frac{\sigma-2(m+1)}{3} \rfloor} \times2^{\lfloor\frac{3(m+1)-\sigma}{2} \rfloor}=2\times6^{\lfloor\frac{\sigma-2(m+1)}{3} \rfloor} \times2^{\lfloor\frac{3(m+1)-\sigma}{2} \rfloor}$. This contradicts to ${\rm per}D\geq 2\times6^{\lfloor\frac{\sigma-2(m+1)}{3} \rfloor} \times2^{\lfloor\frac{3(m+1)-\sigma}{2} \rfloor}$ in (\ref{equ18}).
Set that $s>5$ or $t>5$. By (iii) of Claim \ref{cla2}, we have ${\rm per}D<2\times6^{\lfloor\frac{\sigma-2(m+1)}{3} \rfloor} \times2^{\lfloor\frac{3(m+1)-\sigma}{2} \rfloor}$. This contradicts to ${\rm per}D\geq 2\times6^{\lfloor\frac{\sigma-2(m+1)}{3} \rfloor} \times2^{\lfloor\frac{3(m+1)-\sigma}{2} \rfloor}$ in (\ref{equ18}).

Assume that $k_{2}=k_{3}=4$ and $\sigma-2(m+1)\equiv1(mod~3)$, and $3(m+1)-\sigma$ is odd. It implies that   $\sigma-2(m+1)-s\equiv2(mod~3)$ and $3(m+1)-\sigma+s-1$ is odd. Since $\sigma-2(m+1)-s\equiv2(mod~3)$ and $3(m+1)-\sigma+s-1$ is odd, we have $s=2,5,8,\cdots$, and  $s-1=1,3,5,\cdots$. This means that $s=2,8,14\cdots$. Similarly, we have  $t=2,8,14\cdots$. Set that $s=2$ and $t=2$. By (\ref{equ18}), we have ${\rm per}A={\rm per}D=4\times6^{\lfloor\frac{\sigma-2(m+1)-s}{3} \rfloor}\times2^{\lfloor\frac{3(m+1)-\sigma+s-1}{2} \rfloor}+4\times6^{\lfloor\frac{\sigma-2(m+1)-t}{3} \rfloor} \times2^{\lfloor\frac{3(m+1)-\sigma+t-1}{2} \rfloor}=4\times2\times6^{\lfloor\frac{\sigma-2(m+1)-2}{3} \rfloor}\times2^{\lfloor\frac{3(m+1)-\sigma+1}{2} \rfloor}=4\times2\times4\times\frac{1}{6}\times6^{\lfloor\frac{\sigma-2(m+1)}{3} \rfloor} \times2^{\lfloor\frac{3(m+1)-\sigma}{2} \rfloor}=\frac{4}{3}\times6^{\lfloor\frac{\sigma-2(m+1)}{3} \rfloor} \times2^{\lfloor\frac{3(m+1)-\sigma}{2} \rfloor}<\frac{3}{2}\times6^{\lfloor\frac{\sigma-2(m+1)}{3} \rfloor} \times2^{\lfloor\frac{3(m+1)-\sigma}{2} \rfloor}=\frac{3}{2}\times6^{\lfloor\frac{\sigma-2(m+1)}{3} \rfloor} \times2^{\lfloor\frac{3(m+1)-\sigma}{2} \rfloor}$.  This contradicts to ${\rm per}D\geq \frac{3}{2}\times6^{\lfloor\frac{\sigma-2(m+1)}{3} \rfloor} \times2^{\lfloor\frac{3(m+1)-\sigma}{2} \rfloor}$ in (\ref{equ18}).
Set that $s>2$ or $t>2$.  By (iii) of Claim \ref{cla2}, we have ${\rm per}D<\frac{3}{2}\times6^{\lfloor\frac{\sigma-2(m+1)}{3} \rfloor} \times2^{\lfloor\frac{3(m+1)-\sigma}{2} \rfloor}$, a contradiction.

Assume that $k_{2}=k_{3}=4$ and $\sigma-2(m+1)\equiv0(mod~3)$. It implies that  $\sigma-2(m+1)-s\equiv2(mod~3)$ and $3(m+1)-\sigma+s-1$
is odd.  This means that $s=1,4,7,\cdots$. Similarly, we have $t=1,4,7,\cdots$. Set that $s=1$ and $t=1$. By (\ref{equ18}), we have ${\rm per}A={\rm per}D=4\times6^{\lfloor\frac{\sigma-2(m+1)-s}{3} \rfloor} \times2^{\lfloor\frac{3(m+1)-\sigma+s-1}{2} \rfloor}+4\times6^{\lfloor\frac{\sigma-2(m+1)-t}{3} \rfloor} \times2^{\lfloor\frac{3(m+1)-\sigma+t-1}{2} \rfloor}=4\times2\times6^{\lfloor\frac{\sigma-2(m+1)-1}{3} \rfloor} \times2^{\lfloor\frac{3(m+1)-\sigma}{2} \rfloor}=4\times2\times4\times\frac{1}{6}\times6^{\lfloor\frac{\sigma-2(m+1)}{3} \rfloor} \times2^{\lfloor\frac{3(m+1)-\sigma}{2} \rfloor}=\frac{4}{3}\times6^{\lfloor\frac{\sigma-2(m+1)}{3} \rfloor} \times2^{\lfloor\frac{3(m+1)-\sigma}{2} \rfloor}>6^{\lfloor\frac{\sigma-2(m+1)}{3} \rfloor} \times2^{\lfloor\frac{3(m+1)-\sigma}{2} \rfloor}=6^{\lfloor\frac{\sigma-2(m+1)}{3} \rfloor} \times2^{\lfloor\frac{3(m+1)-\sigma}{2} \rfloor}$. Set  $s=1$ or $t=1$. Since $\sigma-2(m+1)-s=\sigma-2(m+1)-1\equiv2(mod~3)$ and $3(m+1)-\sigma+s-1=3(m+1)-\sigma$ is odd, we get that $m$-square submatrix $D(1,1)$ satisfies item (III) in Theorem \ref{the2}. However,   $m$-square submatrix $D(1,1)$ is not combinatorially equivalent to $N_{m,\sigma}$. This  violates  the induction hypothesis. Set that $s>1$ or $s>1$. By (iv) of Claim \ref{cla2}, we have ${\rm per}D<6^{\lfloor\frac{\sigma-2(m+1)}{3} \rfloor} \times2^{\lfloor\frac{3(m+1)-\sigma}{2} \rfloor}$. This contradicts to ${\rm per}D\geq 6^{\lfloor\frac{\sigma-2(m+1)}{3} \rfloor} \times2^{\lfloor\frac{3(m+1)-\sigma}{2} \rfloor}$ in (\ref{equ18}).

Assume that $k_{2}=k_{3}=3$ and  $\sigma-2(m+1)\equiv2(mod~3)$, and $3(m+1)-\sigma$ is odd. This implies that $\sigma-2(m+1)-s\equiv2(mod~3)$ and $3(m+1)-\sigma+s-1$. Since $\sigma-2(m+1)-s\equiv2(mod~3)$ and $3(m+1)-\sigma+s-1$ is even, we get that $s=3,6,9,\cdots$, and $s-1=1,3,5,\cdots$. This means that $s=6,12,18,\cdots$. Similarly, we have $t=6,12,18,\cdots$.
Set that $s=6$ and $t=6$.  By (\ref{equ18}), we have ${\rm per}A={\rm per}D=3\times6^{\lfloor\frac{\sigma-2(m+1)-s}{3} \rfloor} \times2^{\lfloor\frac{3(m+1)-\sigma+s-1}{2} \rfloor}+3\times6^{\lfloor\frac{\sigma-2(m+1)-t}{3}\rfloor} \times2^{\lfloor\frac{3(m+1)-\sigma+t-1}{2} \rfloor}=\frac{4}{3}\times6^{\lfloor\frac{\sigma-2(m+1)}{3} \rfloor}\times2^{\lfloor\frac{3(m+1)-\sigma}{2} \rfloor}<4\times6^{\lfloor\frac{\sigma-2(m+1)}{3} \rfloor} \times2^{\lfloor\frac{3(m+1)-\sigma}{2} \rfloor}=4\times6^{\lfloor\frac{\sigma-2(m+1)}{3} \rfloor}\times2^{\lfloor\frac{3(m+1)-\sigma}{2} \rfloor}$. This contradicts to ${\rm per}D\geq 4\times6^{\lfloor\frac{\sigma-2(m+1)}{3} \rfloor} \times2^{\lfloor\frac{3(m+1)-\sigma}{2} \rfloor}$ in (\ref{equ18}).
Set that $s>6$ or $t>6$. By (i) of Claim \ref{cla2}, we have ${\rm per}D<4\times6^{\lfloor\frac{\sigma-2(m+1)}{3} \rfloor} \times2^{\lfloor\frac{3(m+1)-\sigma}{2} \rfloor}$, a contradiction.

Assume that $k_{2}=k_{3}=3$ and $\sigma-2(m+1)\equiv2(mod~3)$, and $3(m+1)-\sigma$ is odd.  This implies that   $\sigma-2(m+1)-s\equiv2(mod~3)$ and $3(m+1)-\sigma+s-1$ is even. Since $\sigma-2(m+1)-s\equiv2(mod~3)$ and $3(m+1)-\sigma+s-1$ is even, we have $s=3,6,9,\cdots$, and $s-1=0,2,4,\cdots$. This means that  $s=3,9,15,\cdots$. Similarly, we have $t=3,9,15,\cdots$.  Set $s=3$ and $t=3$.  By (\ref{equ18}), we have ${\rm per}A={\rm per}D=3\times6^{\lfloor\frac{\sigma-2(m+1)-s}{3} \rfloor} \times2^{\lfloor\frac{3(m+1)-\sigma+s-1}{2} \rfloor}+3\times6^{\lfloor\frac{\sigma-2(m+1)-t}{3} \rfloor} \times2^{\lfloor\frac{3(m+1)-\sigma+t-1}{2} \rfloor}=2\times6^{\lfloor\frac{\sigma-2(m+1)}{3} \rfloor} \times2^{\lfloor\frac{3(m+1)-\sigma}{2}\rfloor}<3\times6^{\lfloor\frac{\sigma-2(m+1)}{3} \rfloor}\times2^{\lfloor\frac{3(m+1)-\sigma}{2} \rfloor}=3\times6^{\lfloor\frac{\sigma-2(m+1)}{3} \rfloor} \times2^{\lfloor\frac{3(m+1)-\sigma}{2} \rfloor}$. This contradicts to ${\rm per}D\geq3\times6^{\lfloor\frac{\sigma-2(m+1)}{3} \rfloor} \times2^{\lfloor\frac{3(m+1)-\sigma}{2}\rfloor}$ in (\ref{equ18}). Set that $s>6$ or $t>6$. By (ii) of Claim \ref{cla2}, we have ${\rm per}D<3\times6^{\lfloor\frac{\sigma-2(m+1)}{3} \rfloor} \times2^{\lfloor\frac{3(m+1)-\sigma}{2} \rfloor}$, a contradiction.

Assume that $k_{2}=k_{3}=3$ and $\sigma-2(m+1)\equiv1(mod~3)$, and $3(m+1)-\sigma$ is odd. This implies that   $\sigma-2(m+1)-s\equiv2(mod~3)$, and $3(m+1)-\sigma+s-1$ is even. Since $\sigma-2(m+1)-s\equiv2(mod~3)$ and $3(m+1)-\sigma+s-1$ is even, we have  $s=2,5,8,\cdots$, and $s-1=1,3,5,\cdots$. This means that $s=2,8,14\cdots$.  Similarly, we have $t=2,8,14\cdots$. Set that  $s=2$ and $t=2$. By (\ref{equ18}), we have ${\rm per}A={\rm per}D=3\times6^{\lfloor\frac{\sigma-2(m+1)-s}{3} \rfloor} \times2^{\lfloor\frac{3(m+1)-\sigma+s-1}{2} \rfloor}+3\times6^{\lfloor\frac{\sigma-2(m+1)-t}{3}\rfloor}\times2^{\lfloor\frac{3(m+1)-\sigma+t-1}{2} \rfloor}=2\times6^{\lfloor\frac{\sigma-2(m+1)}{3}\rfloor}\times2^{\lfloor\frac{3(m+1)-\sigma}{2} \rfloor}=2\times6^{\lfloor\frac{\sigma-2(m+1)}{3} \rfloor}\times2^{\lfloor\frac{3(m+1)-\sigma}{2}\rfloor}$. Since $\sigma-2(m+1)-s=\sigma-2(m+1)-2\equiv2(mod~3)$, and $3(m+1)-\sigma+s-1=3(m+1)-\sigma+1$ is even, we get that $m$-square submatrix $D(1,1)$ satisfies item (III) in Theorem \ref{the2}. This implies that $m$-square submatrix $D(1,1)$ is combinatorially equivalent to $M_{m,\sigma}$. Hence,  $(m+1)$-square matrix $A$ is combinatorially equivalent to $T_{(m+1),\sigma}$,
\\
$T^{'}=\left[
\begin{matrix}
\begin{array}{c|cc|c}

1    &  1 &  0 &   0\cdots0 \\
0    &  1 &  1 &   0\cdots0 \\
0    &  1 &  1 &   0\cdots0 \\

\vdots&  \vdots&  \vdots&  \\

1      &   \vdots &  \vdots &     \\

\vdots & \vdots   & \vdots  &    C \\

 1     &  \vdots   & \vdots &       \\

\vdots & \vdots & \vdots &       \\

\end{array}
\end{matrix}\right]$
or $T^{''}=\left[
\begin{matrix}
\begin{array}{c|cc|c}

1    &  1 &  0 &   0\cdots0 \\
0    &  0 &  1 &   0\cdots0 \\
0    &  0 &  1 &   0\cdots0 \\

\vdots&  \vdots&  \vdots&  \\

1      &   1  &  \vdots &     \\

\vdots & \vdots   & \vdots  &    C \\

 1     &  1   & \vdots &       \\

\vdots & \vdots & \vdots &       \\

\end{array}
\end{matrix}\right]$.
If $A=T^{'} {\rm or}~ T^{''}$, then ${\rm per}A\neq\mu((m+1),\tau)= 4\times6^{\lfloor\frac{\sigma-2(m+1)}{3} \rfloor} \times2^{\lfloor\frac{3(m+1)-\sigma}{2} \rfloor}$. Hence,  $(m+1)$-square matrix $A$ is combinatorially equivalent to $N_{(m+1),\sigma}$.
Set that  $s>2$ or $t>2$. By (iii) of Claim \ref{cla2}, we have ${\rm per}D<2\times6^{\lfloor\frac{\sigma-2(m+1)}{3} \rfloor} \times2^{\lfloor\frac{3(m+1)-\sigma}{2} \rfloor}$. This contradicts to ${\rm per}D\geq2\times6^{\lfloor\frac{\sigma-2(m+1)}{3} \rfloor} \times2^{\lfloor\frac{3(m+1)-\sigma}{2}\rfloor}$ in (\ref{equ18}).

Assume that $k_{2}=k_{3}=3$ and $\sigma-2(m+1)\equiv1(mod~3)$, and $3(m+1)-\sigma$ is even. It implies that  $\sigma-2(m+1)-s\equiv2(mod~3)$ and $3(m+1)-\sigma+s-1$ is even. Since $\sigma-2(m+1)-s\equiv2(mod~3)$ and $3(m+1)-\sigma+s-1$ is even, we have  $s=2,5,8,\cdots$, and $s-1=0,2,4,\cdots$. This means that $s=5,11,17\cdots$. Similarly, we have $t=5,11,17\cdots$. Set that $s=5$ and $t=5$. By (\ref{equ18}), we have ${\rm per}A={\rm per}D=3\times6^{\lfloor\frac{\sigma-2(m+1)-s}{3} \rfloor}\times2^{\lfloor\frac{3(m+1)-\sigma+s-1}{2} \rfloor}+3\times6^{\lfloor\frac{\sigma-2(m+1)-t}{3} \rfloor} \times2^{\lfloor\frac{3(m+1)-\sigma+t-1}{2} \rfloor}=\frac{2}{3}\cdot6^{\lfloor\frac{\sigma-2(m+1)}{3} \rfloor} \times2^{\lfloor\frac{3(m+1)-\sigma}{2} \rfloor}<\frac{3}{2}\cdot6^{\lfloor\frac{\sigma-2(m+1)}{3} \rfloor} \times2^{\lfloor\frac{3(m+1)-\sigma}{2} \rfloor}=\frac{3}{2}\times6^{\lfloor\frac{\sigma-2(m+1)}{3} \rfloor} \times2^{\lfloor\frac{3(m+1)-\sigma}{2} \rfloor}$. This contradicts to ${\rm per}D\geq\frac{3}{2}\times6^{\lfloor\frac{\sigma-2(m+1)}{3} \rfloor} \times2^{\lfloor\frac{3(m+1)-\sigma}{2} \rfloor}$ in (\ref{equ18}). Set that $s>5$ or $t>5$. By (iv) of Claim \ref{cla2}, we have ${\rm per}D<\frac{3}{2}\times6^{\lfloor\frac{\sigma-2(m+1)}{3} \rfloor} \times2^{\lfloor\frac{3(m+1)-\sigma}{2} \rfloor}$, a contradiction.

Assume that $k_{2}=k_{3}=3$ and $\sigma-2(m+1)\equiv0(mod~3)$. It implies that  $\sigma-2(m+1)-s\equiv2(mod~3)$ and $3(m+1)-\sigma+s-1$ is even. Since $\sigma-2(m+1)-s\equiv2(mod~3)$ and $\sigma-2(m+1)\equiv0(mod~3)$,  we have $s=1,4,7,\cdots$. Similarly, we have $t=1,4,7,\cdots$. Set that $s=1$ and $t=1$.  By (\ref{equ18}), we have ${\rm per}A={\rm per}D=3\times6^{\lfloor\frac{\sigma-2(m+1)-s}{3} \rfloor} \times2^{\lfloor\frac{3(m+1)-\sigma+s-1}{2} \rfloor}+3\times6^{\lfloor\frac{\sigma-2(m+1)-t}{3} \rfloor} \times2^{\lfloor\frac{3(m+1)-\sigma+t-1}{2} \rfloor}=6^{\lfloor\frac{\sigma-2(m+1)}{3} \rfloor}\times2^{\lfloor\frac{3(m+1)-\sigma}{2} \rfloor}=6^{\lfloor\frac{\sigma-2(m+1)}{3}\rfloor} \times2^{\lfloor\frac{3(m+1)-\sigma}{2} \rfloor}$. Since $\sigma-2(m+1)-s\equiv2(mod~3)$ and $3(m+1)-\sigma+s-1$ is even,  we get that  $m$-square  submatrix $D(1,1)$ satisfies item (III) in Theorem \ref{the2}. This means that  $m$-square submatrix $D(1,1)$ is combinatorially equivalent to $M_{m,\sigma}$. Hence,  $(m+1)$-square matrix $A$ is combinatorially equivalent to $F_{(m+1),\sigma}$.
Set that $s>1$ or $t>1$. By (iii) of Claim \ref{cla2}, we have ${\rm per}D<6^{\lfloor\frac{\sigma-2(m+1)}{3} \rfloor} \times2^{\lfloor\frac{3(m+1)-\sigma}{2} \rfloor}$. This contradicts to ${\rm per}D\geq6^{\lfloor\frac{\sigma-2(m+1)}{3} \rfloor} \times2^{\lfloor\frac{3(m+1)-\sigma}{2} \rfloor}$ in (\ref{equ18}).

Assume that $k_{2}=k_{3}=2$ and $\sigma-2(m+1)\equiv2(mod~3)$, and $3(m+1)-\sigma$ is odd.  This implies that  $\sigma-2(m+1)-s\equiv1(mod~3)$ and $3(m+1)-\sigma+s-1$. Since $\sigma-2(m+1)-s\equiv1(mod~3)$ and $3(m+1)-\sigma+s-1$ is odd, we obtain that $s=1,4,7,\cdots$, and $s-1=0,2,4,\cdots$. This implies that $s=1,7,13,\cdots$. Similarly, we have $t=1,7,13,\cdots$. Set that $s=1$ and $t=1$. By (\ref{equ18}), we have ${\rm per}A={\rm per}D=2\times6^{\lfloor\frac{\sigma-2(m+1)-s}{3}\rfloor}\times2^{\lfloor\frac{3(m+1)-\sigma+s-1}{2} \rfloor}+2\times6^{\lfloor\frac{\sigma-2(m+1)-t}{3} \rfloor}\times2^{\lfloor\frac{3(m+1)-\sigma+t-1}{2} \rfloor}=2\times2\times6^{\lfloor\frac{\sigma-2(m+1)-1}{3}\rfloor}\times2^{\lfloor\frac{3(m+1)-\sigma}{2} \rfloor}=4\times6^{\lfloor\frac{\sigma-2(m+1)}{3}\rfloor}\times2^{\lfloor\frac{3(m+1)-\sigma}{2} \rfloor}=4\times6^{\lfloor\frac{\sigma-2(m+1)}{3} \rfloor}\times2^{\lfloor\frac{3(m+1)-\sigma}{2}\rfloor}$.  Since $s=1$, $\sigma-2(m+1)-s=\sigma-2(m+1)-1\equiv1(mod~3)$ and $3(m+1)-\sigma+s-1=3(m+1)-\sigma$ is odd,   we get that  $m$-square  submatrix $D(1,1)$ satisfies item (II) in Theorem \ref{the2}. This implies taht $m$-square submatrix $D(1,1)$ is combinatorially equivalent to $T_{m,\sigma}$. Furthermore,  $(m+1)$-square matrix $A$ is combinatorially equivalent to $N_{(m+1),\sigma}$.  Set that $s=7$ and $t=1$. By (\ref{equ18}), we have ${\rm per}A={\rm per}D=2\times6^{\lfloor\frac{\sigma-2(m+1)-s}{3} \rfloor}\times2^{\lfloor\frac{3(m+1)-\sigma+s-1}{2} \rfloor}+2\cdot6^{\lfloor\frac{\sigma-2(m+1)-t}{3} \rfloor} \times2^{\lfloor\frac{3(m+1)-\sigma+t-1}{2} \rfloor}=2\times6^{\lfloor\frac{\sigma-2(m+1)-7}{3} \rfloor} \times2^{\lfloor\frac{3(m+1)-\sigma+6}{2} \rfloor}+2\times6^{\lfloor\frac{\sigma-2(m+1)}{3} \rfloor} \times2^{\lfloor\frac{3(m+1)-\sigma}{2} \rfloor}=2\times\frac{1}{36}\times8\times6^{\lfloor\frac{\sigma-2(m+1)}{3} \rfloor} \times2^{\lfloor\frac{3(m+1)-\sigma}{2} \rfloor}+2\times6^{\lfloor\frac{\sigma-2(m+1)}{3} \rfloor} \times2^{\lfloor\frac{3(m+1)-\sigma}{2} \rfloor}=\frac{22}{9}\times6^{\lfloor\frac{\sigma-2(m+1)}{3} \rfloor} \times2^{\lfloor\frac{3(m+1)-\sigma}{2} \rfloor}<4\times6^{\lfloor\frac{\sigma-2(m+1)}{3} \rfloor} \times2^{\lfloor\frac{3(m+1)-\sigma}{2} \rfloor}$.  This contradicts to ${\rm per}D\geq 4\times6^{\lfloor\frac{\sigma-2(m+1)}{3} \rfloor} \times2^{\lfloor\frac{3(m+1)-\sigma}{2} \rfloor}$ in (\ref{equ18}).
Set that $s>7$ or $t>1$. By (i) of Claim \ref{cla2}, we have ${\rm per}D<4\times6^{\lfloor\frac{\sigma-2(m+1)}{3} \rfloor} \times2^{\lfloor\frac{3(m+1)-\sigma}{2} \rfloor}$, a contradiction.

Assume that $k_{2}=k_{3}=2$ and $\sigma-2(m+1)\equiv2(mod~3)$, and $3(m+1)-\sigma$ is even. This implies that  $\sigma-2(m+1)-s\equiv1(mod~3)$ and $3(m+1)-\sigma+s-1$ is odd. Since $\sigma-2(m+1)-s\equiv1(mod~3)$ and $3(m+1)-\sigma+s-1$ is odd, we have  $s=1,4,7,\cdots$, and $s-1=1,3,5,\cdots$. This means that $s=4,10,16,\cdots$. Similarly, we have $t=4,10,16,\cdots$. Set that $s=4$ and $t=4$. By (\ref{equ18}), we have ${\rm per}A={\rm per}D=2\times6^{\lfloor\frac{\sigma-2(m+1)-s}{3} \rfloor} \times2^{\lfloor\frac{3(m+1)-\sigma+s-1}{2} \rfloor}+2\times6^{\lfloor\frac{\sigma-2(m+1)-t}{3} \rfloor}\times2^{\lfloor\frac{3(m+1)-\sigma+t-1}{2} \rfloor}=\frac{4}{3}\times6^{\lfloor\frac{\sigma-2(m+1)}{3}\rfloor} \times2^{\lfloor\frac{3(m+1)-\sigma}{2} \rfloor}<2\times6^{\lfloor\frac{\sigma-2(m+1)}{3} \rfloor}\times2^{\lfloor\frac{3(m+1)-\sigma}{2} \rfloor}=3\times6^{\lfloor\frac{\sigma-2(m+1)}{3}\rfloor}\times2^{\lfloor\frac{3(m+1)-\sigma}{2} \rfloor}$. This contradicts to ${\rm per}D\geq 3\times6^{\lfloor\frac{\sigma-2(m+1)}{3}\rfloor}\times2^{\lfloor\frac{3(m+1)-\sigma}{2} \rfloor}$ in (\ref{equ18}).
Set that $s>4$ or $t>4$. By (ii) of Claim \ref{cla2}, we have ${\rm per}D<3\times6^{\lfloor\frac{\sigma-2(m+1)}{3} \rfloor} \times2^{\lfloor\frac{3(m+1)-\sigma}{2} \rfloor}$, a contradiction.

Assume that $k_{2}=k_{3}=2$ and $\sigma-2(m+1)\equiv1(mod~3)$, and $3(m+1)-\sigma$ is odd. It implies that  $\sigma-2(m+1)-s\equiv1(mod~3)$ and $3(m+1)-\sigma+s-1$ is odd. Since $\sigma-2(m+1)-s\equiv1(mod~3)$ and $3(m+1)-\sigma+s-1$ is odd, we have $s=3,6,9,\cdots$, and $s-1=0,2,4,\cdots$. This means that $s=3,9,15,\cdots$. Similarly, we have $t=3,9,15,\cdots$.
Set that $s=3$ and $t=3$. By (\ref{equ18}), we have ${\rm per}A={\rm per}D=2\times6^{\lfloor\frac{\sigma-2(m+1)-s}{3} \rfloor} \times2^{\lfloor\frac{3(m+1)-\sigma+s-1}{2} \rfloor}+2\times6^{\lfloor\frac{\sigma-2(m+1)-t}{3} \rfloor} \times2^{\lfloor\frac{3(m+1)-\sigma+t-1}{2} \rfloor}=\frac{4}{3}\times6^{\lfloor\frac{\sigma-2(m+1)}{3} \rfloor} \times2^{\lfloor\frac{3(m+1)-\sigma}{2} \rfloor}<2\times6^{\lfloor\frac{\sigma-2(m+1)}{3} \rfloor} \times2^{\lfloor\frac{3(m+1)-\sigma}{2} \rfloor}=2\times6^{\lfloor\frac{\sigma-2(m+1)}{3} \rfloor} \times2^{\lfloor\frac{3(m+1)-\sigma}{2} \rfloor}$. This contradicts to ${\rm per}C\geq 2\times6^{\lfloor\frac{\sigma-2(m+1)}{3} \rfloor} \times2^{\lfloor\frac{3(m+1)-\sigma}{2} \rfloor}$ in (\ref{equ18}).
Set $s>3$ or $t>3$. By (iii) of Claim \ref{cla2}, we have ${\rm per}D<2\times6^{\lfloor\frac{\sigma-2(m+1)}{3} \rfloor} \times2^{\lfloor\frac{3(m+1)-\sigma}{2} \rfloor}$, a contradiction.

Assume that $k_{2}=k_{3}=2$ and $\sigma-2(m+1)\equiv1(mod~3)$, and $3(m+1)-\sigma$ is even. It implies that $\sigma-2(m+1)-s\equiv1(mod~3)$ and $3(m+1)-\sigma+s-1$ is odd.  Since $\sigma-2(m+1)-s\equiv1(mod~3)$ and $3(m+1)-\sigma+s-1$ is odd, we have  $s=3,6,9,\cdots$, and $s-1=1,3,5,\cdots$.  This means that  $s=6,12,18\cdots$. Similarly, we have $t=6,12,18\cdots$. Set that $s=6$ and $t=6$.  By (\ref{equ18}), we have ${\rm per}A={\rm per}D=2\times6^{\lfloor\frac{\sigma-2(m+1)-s}{3} \rfloor} \times2^{\lfloor\frac{3(m+1)-\sigma+s-1}{2} \rfloor}+2\times6^{\lfloor\frac{\sigma-2(m+1)-t}{3} \rfloor} \times2^{\lfloor\frac{3(m+1)-\sigma+t-1}{2} \rfloor}=\frac{4}{9}\times6^{\lfloor\frac{\sigma-2(m+1)}{3} \rfloor} \times2^{\lfloor\frac{3(m+1)-\sigma}{2} \rfloor}<\frac{3}{2}\times6^{\lfloor\frac{\sigma-2(m+1)}{3} \rfloor} \times2^{\lfloor\frac{3(m+1)-\sigma}{2} \rfloor}=\frac{3}{2}\times6^{\lfloor\frac{\sigma-2(m+1)}{3} \rfloor} \times2^{\lfloor\frac{3(m+1)-\sigma}{2} \rfloor}$. This contradicts to ${\rm per}D\geq \frac{3}{2}\times6^{\lfloor\frac{\sigma-2(m+1)}{3} \rfloor} \times2^{\lfloor\frac{3(m+1)-\sigma}{2} \rfloor}$ in (\ref{equ18}). Set that $s>6$ or $t>6$. By (iv) of Claim \ref{cla2}, we have ${\rm per}D<\frac{3}{2}\times6^{\lfloor\frac{\sigma-2(m+1)}{3} \rfloor} \times2^{\lfloor\frac{3(m+1)-\sigma}{2} \rfloor}$, a contradiction.

Assume that $k_{2}=k_{3}=2$ and  $\sigma-2(m+1)\equiv1(mod~3)$. It implies that $3(m+1)-\sigma$ is even,  $\sigma-2(m+1)-s\equiv1(mod~3)$ and $3(m+1)-\sigma+s-1$ is odd. Since $\sigma-2(m+1)-s\equiv1(mod~3)$ and $3(m+1)-\sigma+s-1$ is odd,  we have $s=2,5,8,\cdots$. Similarly, we have $t=2,5,8,\cdots$. Set that $s=2$ and $t=2$. By (\ref{equ18}), we have ${\rm per}A={\rm per}D=2\times6^{\lfloor\frac{\sigma-2(m+1)-s}{3} \rfloor}\times2^{\lfloor\frac{3(m+1)-\sigma+s-1}{2} \rfloor}+2\times6^{\lfloor\frac{\sigma-2(m+1)-t}{3} \rfloor}\times2^{\lfloor\frac{3(m+1)-\sigma+t-1}{2} \rfloor}=\frac{2}{3}\times6^{\lfloor\frac{\sigma-2(m+1)}{3}\rfloor} \times2^{\lfloor\frac{3(m+1)-\sigma}{2} \rfloor}<6^{\lfloor\frac{\sigma-2(m+1)}{3} \rfloor}\times2^{\lfloor\frac{3(m+1)-\sigma}{2} \rfloor}=6^{\lfloor\frac{\sigma-2(m+1)}{3} \rfloor} \times2^{\lfloor\frac{3(m+1)-\sigma}{2} \rfloor}$. This contradicts to ${\rm per}D\geq 6^{\lfloor\frac{\sigma-2(m+1)}{3} \rfloor} \times2^{\lfloor\frac{3(m+1)-\sigma}{2} \rfloor}$ in (\ref{equ18}).
Set that  $s>2$ or  $s>2$. By (iii) of Claim \ref{cla2}, we have ${\rm per}D<6^{\lfloor\frac{\sigma-2(m+1)}{3} \rfloor} \times2^{\lfloor\frac{3(m+1)-\sigma}{2} \rfloor}$, a contradiction.

Assume that $k_{2}=k_{3}=\frac{3}{2}$ and $\sigma-2(m+1)\equiv2(mod~3)$, and $3(m+1)-\sigma$ is odd. This implies that $\sigma-2(m+1)-s\equiv1(mod~3)$ and $3(m+1)-\sigma+s-1$ is even. Since $\sigma-2(m+1)-s\equiv1(mod~3)$ and $3(m+1)-\sigma+s-1$ is even, we have $s=2,5,8,\cdots$, and $s-1=1,3,5,\cdots$.  This means that $s=2,8,14,\cdots$. Similarly, we have $t=2,8,14,\cdots$. Set that $s=2$ and $t=2$. By (\ref{equ18}), we have ${\rm per}A={\rm per}D=\frac{3}{2}\times6^{\lfloor\frac{\sigma-2(m+1)-s}{3} \rfloor} \times2^{\lfloor\frac{3(m+1)-\sigma+s-1}{2} \rfloor}+\frac{3}{2}\times6^{\lfloor\frac{\sigma-2(m+1)-t}{3} \rfloor} \times2^{\lfloor\frac{3(m+1)-\sigma+t-1}{2} \rfloor}=2\times\frac{3}{2}\times6^{\lfloor\frac{\sigma-2(m+1)-2}{3} \rfloor} \times2^{\lfloor\frac{3(m+1)-\sigma+1}{2} \rfloor}=3\times6^{\lfloor\frac{\sigma-2(m+1)}{3} \rfloor} \times2^{\lfloor\frac{3(m+1)-\sigma}{2} \rfloor}<4\times6^{\lfloor\frac{\sigma-2(m+1)}{3} \rfloor} \times2^{\lfloor\frac{3(m+1)-\sigma}{2} \rfloor}$.  This contradicts to ${\rm per}D\geq 4\times6^{\lfloor\frac{\sigma-2(m+1)}{3} \rfloor} \times2^{\lfloor\frac{3(m+1)-\sigma}{2} \rfloor}$ in (\ref{equ18}).
Set that  $s>2$ or  $s>2$. By (iii) of Claim \ref{cla2}, we have ${\rm per}D<4\times6^{\lfloor\frac{\sigma-2(m+1)}{3} \rfloor} \times2^{\lfloor\frac{3(m+1)-\sigma}{2} \rfloor}$, a contradiction.

Assume that $k_{2}=k_{3}=\frac{3}{2}$ and $\sigma-2(m+1)\equiv2(mod~3)$, and $3(m+1)-\sigma$ is even. This implies that $\sigma-2(m+1)-s\equiv1(mod~3)$ and $3(m+1)-\sigma+s-1$ is even. Since $\sigma-2(m+1)-s\equiv1(mod~3)$ and $3(m+1)-\sigma+s-1$ is even, we have $s=1,4,7,\cdots$, and $s-1=0,2,4,\cdots$.  This means that $s=1,7,13,\cdots$. Similarly, we have $t=1,7,13,\cdots$. Set that $s=1$ and $t=1$. By (\ref{equ18}), we have ${\rm per}A={\rm per}D=\frac{3}{2}\times6^{\lfloor\frac{\sigma-2(m+1)-s}{3} \rfloor} \times2^{\lfloor\frac{3(m+1)-\sigma+s-1}{2} \rfloor}+\frac{3}{2}\times6^{\lfloor\frac{\sigma-2(m+1)-t}{3} \rfloor} \times2^{\lfloor\frac{3(m+1)-\sigma+t-1}{2} \rfloor}=3\times6^{\lfloor\frac{\sigma-2(m+1)}{3} \rfloor} \times2^{\lfloor\frac{3(m+1)-\sigma}{2} \rfloor}=3\times6^{\lfloor\frac{\sigma-2(m+1)}{3}\rfloor}\times2^{\lfloor\frac{3(m+1)-\sigma}{2} \rfloor}$. Since $\sigma-2(m+1)-s=\sigma-2(m+1)-2\equiv0(mod~3)$ and $3(m+1)-\sigma+s-1=3n-\sigma+1$ is even, we get that  $m$-square  submatrix $D(1,1)$ satisfies item (II) in Theorem \ref{the2}. This implies that $m$-square  submatrix $D(1,1)$ is combinatorially equivalent to $S_{m,\sigma}$. Hence, $(m+1)$-square  matrix $A$ is combinatorially equivalent to $M_{(m+1),\sigma}$. Set that $s>1$ or $t>1$. By (ii) of Claim \ref{cla2}, we have ${\rm per}D<3\times6^{\lfloor\frac{\sigma-2(m+1)}{3} \rfloor} \times2^{\lfloor\frac{3(m+1)-\sigma}{2} \rfloor}$, a contradiction.

Assume that $k_{2}=k_{3}=\frac{3}{2}$ and  $\sigma-2(m+1)\equiv1(mod~3)$, and $3(m+1)-\sigma$ is odd. It implies that   $\sigma-2(m+1)-s\equiv1(mod~3)$ and $3(m+1)-\sigma+s-1$ is even. Since $\sigma-2(m+1)-s\equiv1(mod~3)$ and $3(m+1)-\sigma+s-1$ is even, we have $s=3,6,9,\cdots$, and $s-1=1,3,5,\cdots$. This means that $s=6,12,18\cdots$. Similarly, we have $t=6,12,18\cdots$.
Set that $s=6$ and $t=6$. By (\ref{equ18}), we have ${\rm per}A={\rm per}D=\frac{3}{2}\times6^{\lfloor\frac{\sigma-2(m+1)-s}{3} \rfloor} \times2^{\lfloor\frac{3(m+1)-\sigma+s-1}{2} \rfloor}+\frac{3}{2}\times6^{\lfloor\frac{\sigma-2(m+1)-t}{3} \rfloor} \times2^{\lfloor\frac{3(m+1)-\sigma+t-1}{2} \rfloor}=\frac{2}{3}\times6^{\lfloor\frac{\sigma-2(m+1)}{3} \rfloor} \times2^{\lfloor\frac{3(m+1)-\sigma}{2} \rfloor}<2\times6^{\lfloor\frac{\sigma-2(m+1)}{3} \rfloor} \times2^{\lfloor\frac{3(m+1)-\sigma}{2} \rfloor}=2\times6^{\lfloor\frac{\sigma-2(m+1)}{3} \rfloor} \times2^{\lfloor\frac{3(m+1)-\sigma}{2} \rfloor}$.  This contradicts to ${\rm per}C\geq 2\times6^{\lfloor\frac{\sigma-2(m+1)}{3} \rfloor} \times2^{\lfloor\frac{3(m+1)-\sigma}{2} \rfloor}$ in (\ref{equ18}).
Set that $s>6$ or $t>6$. By (iii) of Claim \ref{cla2}, we have ${\rm per}D<2\times6^{\lfloor\frac{\sigma-2(m+1)}{3} \rfloor} \times2^{\lfloor\frac{3(m+1)-\sigma}{2} \rfloor}$, a contradiction.

Assume that $k_{2}=k_{3}=\frac{3}{2}$ and  $\sigma-2(m+1)\equiv1(mod~3)$, and  $3(m+1)-\sigma$ is even. It implies that  $\sigma-2(m+1)-s\equiv1(mod~3)$ and $3(m+1)-\sigma+s-1$ is even.  Since $\sigma-2(m+1)-s\equiv1(mod~3)$ and $3(m+1)-\sigma+s-1$ is even, we have $s=3,6,9,\cdots$, and $s-1=0,2,4,\cdots$. This means that $s=3,9,15\cdots$. Similarly, we have $t=3,9,15\cdots$. Set that $s=3$ and $t=3$. By (\ref{equ18}), we have ${\rm per}A={\rm per}D=\frac{3}{2}\times6^{\lfloor\frac{\sigma-2(m+1)-s}{3}\rfloor} \times2^{\lfloor\frac{3(m+1)-\sigma+s-1}{2}\rfloor}+\frac{3}{2}\times6^{\lfloor\frac{\sigma-2(m+1)-t}{3}\rfloor} \times2^{\lfloor\frac{3(m+1)-\sigma+t-1}{2}\rfloor}=6^{\lfloor\frac{\sigma-2(m+1)}{3} \rfloor} \times2^{\lfloor\frac{3(m+1)-\sigma}{2} \rfloor}<\frac{3}{2}\times6^{\lfloor\frac{\sigma-2(m+1)}{3}\rfloor} \times2^{\lfloor\frac{3(m+1)-\sigma}{2} \rfloor}=\frac{3}{2}\times6^{\lfloor\frac{\sigma-2(m+1)}{3} \rfloor} \times2^{\lfloor\frac{3(m+1)-\sigma}{2} \rfloor}$. This contradicts to ${\rm per}D\geq \frac{3}{2}\times6^{\lfloor\frac{\sigma-2(m+1)}{3} \rfloor} \times2^{\lfloor\frac{3(m+1)-\sigma}{2} \rfloor}$ in (\ref{equ18}).
Set that $s>3$ or $t>3$. By (iv) of Claim \ref{cla2}, we have ${\rm per}D<\frac{3}{2}\times6^{\lfloor\frac{\sigma-2(m+1)}{3} \rfloor} \times2^{\lfloor\frac{3(m+1)-\sigma}{2} \rfloor}$, a contradiction.

Assume that $k_{2}=k_{3}=\frac{3}{2}$ and $\sigma-2(m+1)\equiv0(mod~3)$. It implies that  $\sigma-2(m+1)-s\equiv1(mod~3)$. This means that $s=2,5,8,\cdots$. Similarly, we have $t=2,5,8,\cdots$. Set that $s=2$ and $t=2$. Then This contradicts to ${\rm per}D\geq k_{1}\times6^{\lfloor\frac{\sigma-2(m+1)}{3} \rfloor} \times2^{\lfloor\frac{3(m+1)-\sigma}{2} \rfloor}$ by the case when $k_{2}=k_{3}=2$ and $k_{1}=1$. Set that $s>2$ or $s>2$.  By (v) of Claim \ref{cla2}, we have ${\rm per}D<k_{1}\times6^{\lfloor\frac{\sigma-2(m+1)}{3} \rfloor} \times2^{\lfloor\frac{3(m+1)-\sigma}{2} \rfloor}$, a contradiction.

Assume that $k_{2}=k_{3}=1$ and $\sigma-2(m+1)\equiv2(mod~3)$, and $3(m+1)-\sigma$ is odd. This implies that  $\sigma-2(m+1)-s\equiv0(mod~3)$. Since $\sigma-2(m+1)\equiv2(mod~3)$ and $\sigma-2(m+1)-s\equiv0(mod~3)$, we have $s=2,5,8,\cdots$. Similarly, we have  $t=2,5,8,\cdots$. Set $s=2$ and  $t=2$. By (\ref{equ18}), we have ${\rm per}A={\rm per}D=6^{\lfloor\frac{\sigma-2(m+1)-s}{3} \rfloor} \times2^{\lfloor\frac{3(m+1)-\sigma+s-1}{2}\rfloor}+6^{\lfloor\frac{\sigma-2(m+1)-t}{3} \rfloor} \times2^{\lfloor\frac{3(m+1)-\sigma+t-1}{2} \rfloor}=2\times6^{\lfloor\frac{\sigma-2(m+1)-2}{3}\rfloor}\times2^{\lfloor\frac{3(m+1)-\sigma+1}{2} \rfloor}=4\times6^{\lfloor\frac{\sigma-2(m+1)}{3} \rfloor}\times2^{\lfloor\frac{3(m+1)-\sigma}{2} \rfloor}=4\times6^{\lfloor\frac{\sigma-2(m+1)}{3} \rfloor}\times2^{\lfloor\frac{3(m+1)-\sigma}{2} \rfloor}$. Since $\sigma-2(m+1)-s=\sigma-2(m+1)-2\equiv0(mod~3)$ and $3(m+1)-\sigma+s-1=3(m+1)-\sigma+1$ is even,  we get that $m$-square   submatrix $D(1,1)$ satisfies item (I) in Theorem \ref{the2}. Hence, $m$-square  submatrix $D(1,1)$ is combinatorially equivalent to $F_{m,\sigma}$.  Moreover, $(m+1)$-square  matrix $A$ is combinatorially equivalent to $N_{(m+1),\sigma}$ or

 $$N^{'}=\left[
\begin{matrix}
\begin{array}{c|cc|c}

1 &  1 &  0 &   0\cdots0 \\
0 &  1 &  1 &   0\cdots0 \\
0 &  1 &  1 &   0\cdots0 \\
\hline
1      &    0   &  0     &         \\

\vdots& \vdots& \vdots&    C      \\

 1     &    0   &  0     &        \\

\end{array}
\end{matrix}\right].$$
Since ${\rm per}N^{'}\neq\mu((m+1),\tau)= 4\times6^{\lfloor\frac{\sigma-2(m+1)}{3} \rfloor} \times2^{\lfloor\frac{3(m+1)-\sigma}{2} \rfloor}$, we obtain that $(m+1)$-square  matrix $A$ is combinatorially equivalent to $N_{(m+1),\sigma}$. Set that $s=5$ and $t=2$, By (\ref{equ18}), we have ${\rm per}A={\rm per}D=6^{\lfloor\frac{\sigma-2(m+1)-s}{3} \rfloor}\times2^{\lfloor\frac{3(m+1)-\sigma+s-1}{2} \rfloor}+6^{\lfloor\frac{\sigma-2(m+1)-t}{3} \rfloor} \times2^{\lfloor\frac{3(m+1)-\sigma+t-1}{2} \rfloor}=6^{\lfloor\frac{\sigma-2(m+1)-5}{3} \rfloor} \times2^{\lfloor\frac{3(m+1)-\sigma+4}{2} \rfloor}+6^{\lfloor\frac{\sigma-2(m+1)-2}{3} \rfloor} \times2^{\lfloor\frac{3(m+1)-\sigma+1}{2} \rfloor}=\frac{1}{6}\times4\times6^{\lfloor\frac{\sigma-2(m+1)}{3} \rfloor} \times2^{\lfloor\frac{3(m+1)-\sigma}{2} \rfloor}+2\times6^{\lfloor\frac{\sigma-2(m+1)}{3} \rfloor} \times2^{\lfloor\frac{3(m+1)-\sigma}{2} \rfloor}=\frac{8}{3}\times6^{\lfloor\frac{\sigma-2(m+1)}{3} \rfloor} \times2^{\lfloor\frac{3(m+1)-\sigma}{2} \rfloor}<4\times6^{\lfloor\frac{\sigma-2(m+1)}{3} \rfloor} \times2^{\lfloor\frac{3(m+1)-\sigma}{2} \rfloor}$. This contradicts to ${\rm per}D\geq 4\times6^{\lfloor\frac{\sigma-2(m+1)}{3} \rfloor} \times2^{\lfloor\frac{3(m+1)-\sigma}{2} \rfloor}$ in (\ref{equ18}).
Set $s>5$.  By (i) of Claim \ref{cla2}, we have ${\rm per}D<4\times6^{\lfloor\frac{\sigma-2(m+1)}{3} \rfloor} \times2^{\lfloor\frac{3(m+1)-\sigma}{2} \rfloor}$, a contradiction.

Assume that $k_{2}=k_{3}=1$, $\sigma-2(m+1)\equiv2(mod~3)$ and $3(m+1)-\sigma$ is even.  This implies that  $\sigma-2(m+1)-s\equiv0(mod~3)$. This means that $s=2,5,8,\cdots$. Similarly, we have $t=2,5,8,\cdots$. Set that $s=2$ and $t=2$. By (\ref{equ18}), we have  ${\rm per}A={\rm per}D=6^{\lfloor\frac{\sigma-2(m+1)-s}{3} \rfloor} \times2^{\lfloor\frac{3(m+1)-\sigma+s-1}{2} \rfloor}+6^{\lfloor\frac{\sigma-2(m+1)-t}{3} \rfloor} \times2^{\lfloor\frac{3(m+1)-\sigma+t-1}{2} \rfloor}=2\times6^{\lfloor\frac{\sigma-2(m+1)}{3} \rfloor} \times2^{\lfloor\frac{3(m+1)-\sigma}{2} \rfloor}<3\times6^{\lfloor\frac{\sigma-2(m+1)}{3} \rfloor} \times2^{\lfloor\frac{3(m+1)-\sigma}{2} \rfloor}=3\times6^{\lfloor\frac{\sigma-2(m+1)}{3} \rfloor} \times2^{\lfloor\frac{3(m+1)-\sigma}{2} \rfloor}$. This contradicts to ${\rm per}D\geq 3\times6^{\lfloor\frac{\sigma-2(m+1)}{3} \rfloor} \times2^{\lfloor\frac{3(m+1)-\sigma}{2} \rfloor}$ in (\ref{equ18}).
Set that $s>2$ or $t>2$. By (ii) of Claim \ref{cla2}, we have ${\rm per}D<3\times6^{\lfloor\frac{\sigma-2(m+1)}{3} \rfloor} \times2^{\lfloor\frac{3(m+1)-\sigma}{2} \rfloor}$, a contradiction.

Assume that $k_{2}=k_{3}=1$, $\sigma-2(m+1)\equiv1(mod~3)$ and $3(m+1)-\sigma$ is odd. It implies that  $\sigma-2(m+1)-s=\sigma-2(m+1)-1\equiv0(mod~3)$ and $3(m+1)-\sigma+s-1$ is odd. This means that $s=1,4,7,\cdots$. Similarly, we have $t=1,4,7,\cdots$.
Set that $s=1$ and  $t=1$. By (\ref{equ18}), we have ${\rm per}A={\rm per}D=6^{\lfloor\frac{\sigma-2(m+1)-s}{3} \rfloor} \times2^{\lfloor\frac{3(m+1)-\sigma+s-1}{2} \rfloor}+6^{\lfloor\frac{\sigma-2(m+1)-t}{3} \rfloor} \times2^{\lfloor\frac{3(m+1)-\sigma+t-1}{2} \rfloor}=2\times6^{\lfloor\frac{\sigma-2(m+1)-1}{3}\rfloor}\times2^{\lfloor\frac{3(m+1)-\sigma}{2} \rfloor}=2\times6^{\lfloor\frac{\sigma-2(m+1)}{3} \rfloor} \times2^{\lfloor\frac{3(m+1)-\sigma}{2}\rfloor}=2\times6^{\lfloor\frac{\sigma-2(m+1)}{3} \rfloor}\times2^{\lfloor\frac{3(m+1)-\sigma}{2} \rfloor}$. Since $\sigma-2(m+1)-s=\sigma-2(m+1)-1\equiv0(mod~3)$ and $3(m+1)-\sigma+s-1$ is odd, we get that $m$-square  submatrix $D(1,1)$ satisfies item (I) in Theorem \ref{the2}. This means that $m$-square  submatrix $D(1,1)$ is combinatorially equivalent to $F_{m,\sigma}^{*}$. Hence, $(m+1)$-square  matrix $A$ is combinatorially equivalent to $P_{(m+1),\sigma}^{*}$. Set that  $s>1$ or $t>1$.  By (iii) of Claim \ref{cla2}, we have ${\rm per}D<2\times6^{\lfloor\frac{\sigma-2(m+1)}{3} \rfloor}\times2^{\lfloor\frac{3(m+1)-\sigma}{2} \rfloor}$, a contradiction.

Assume that $k_{2}=k_{3}=1$, $\sigma-2(m+1)\equiv1(mod~3)$ and $3(m+1)-\sigma$ is even. It implies that  $\sigma-2(m+1)-s=\sigma-2(m+1)-1\equiv0(mod~3)$ and $3(m+1)-\sigma+s-1=3n-\sigma$ is odd.  This means that $s=1,4,7,\cdots$. Similarly, we have $t=1,4,7,\cdots$.
Set that $s=1$ and $t=1$. By (\ref{equ18}), we have ${\rm per}A={\rm per}D=6^{\lfloor\frac{\sigma-2(m+1)-s}{3} \rfloor} \times2^{\lfloor\frac{3(m+1)-\sigma+s-1}{2} \rfloor}+6^{\lfloor\frac{\sigma-2(m+1)-t}{3} \rfloor} \times2^{\lfloor\frac{3(m+1)-\sigma+t-1}{2} \rfloor}=2\times6^{\lfloor\frac{\sigma-2(m+1)}{3} \rfloor}\times2^{\lfloor\frac{3(m+1)-\sigma}{2} \rfloor}>\frac{3}{2}\times6^{\lfloor\frac{\sigma-2(m+1)}{3} \rfloor}\times2^{\lfloor\frac{3(m+1)-\sigma}{2} \rfloor}=\frac{3}{2}\times6^{\lfloor\frac{\sigma-2(m+1)}{3} \rfloor}\times2^{\lfloor\frac{3(m+1)-\sigma}{2}\rfloor}$. Since $\sigma-2(m+1)-s=\sigma-2(m+1)-1\equiv0(mod~3)$ and $3(m+1)-\sigma+s-1=3(m+1)-\sigma$ is odd, we get that submatrix $D(1,1)$ satisfies item (I) in Theorem \ref{the2}. But $m$-square  submatrix $D(1,1)$ is not combinatorially equivalent to $F_{(m+1),\sigma}$. This  violates  the induction hypothesis.
Set that $s=4$ and $t=1$. By (\ref{equ18}), we have ${\rm per}A={\rm per}D=6^{\lfloor\frac{\sigma-2(m+1)-s}{3} \rfloor} \times2^{\lfloor\frac{3(m+1)-\sigma+s-1}{2} \rfloor}+6^{\lfloor\frac{\sigma-2(m+1)-t}{3} \rfloor} \times2^{\lfloor\frac{3(m+1)-\sigma+t-1}{2} \rfloor}=\frac{4}{3}^{\lfloor\frac{\sigma-2(m+1)}{3} \rfloor}\times2^{\lfloor\frac{3(m+1)-\sigma}{2} \rfloor}<\frac{3}{2}\times6^{\lfloor\frac{\sigma-2(m+1)}{3} \rfloor} \times2^{\lfloor\frac{3(m+1)-\sigma}{2} \rfloor}=\frac{3}{2}\times6^{\lfloor\frac{\sigma-2(m+1)}{3} \rfloor} \times2^{\lfloor\frac{3(m+1)-\sigma}{2} \rfloor}$. This contradicts to ${\rm per}D\geq\frac{3}{2}\times6^{\lfloor\frac{\sigma-2(m+1)}{3} \rfloor} \times2^{\lfloor\frac{3(m+1)-\sigma}{2} \rfloor}$ in (\ref{equ18}). Set that $s>4$ or $s>4$. By (iv) of Claim \ref{cla2}, we have ${\rm per}D<\frac{3}{2}\times6^{\lfloor\frac{\sigma-2(m+1)}{3} \rfloor} \times2^{\lfloor\frac{3(m+1)-\sigma}{2} \rfloor}$, a contradiction.

Assume that $k_{2}=k_{3}=1$ and $\sigma-2(m+1)\equiv0(mod~3)$. It implies that $\sigma-2(m+1)-s\equiv0(mod~3)$. This means that $s=3,6,9,\cdots$. Similarly, we have $t=3,6,9,\cdots$. Set that $s=3$ and $t=3$. By (\ref{equ18}), we have ${\rm per}A={\rm per}D=6^{\lfloor\frac{\sigma-2(m+1)-s}{3} \rfloor}\times2^{\lfloor\frac{3(m+1)-\sigma+s-1}{2} \rfloor}+6^{\lfloor\frac{\sigma-2(m+1)-t}{3} \rfloor}\times2^{\lfloor\frac{3(m+1)-\sigma+t-1}{2} \rfloor}=\frac{2}{3}\times6^{\lfloor\frac{\sigma-2(m+1)}{3} \rfloor}\times2^{\lfloor\frac{3(m+1)-\sigma}{2} \rfloor}<6^{\lfloor\frac{\sigma-2(m+1)}{3} \rfloor} \times2^{\lfloor\frac{3(m+1)-\sigma}{2} \rfloor}=6^{\lfloor\frac{\sigma-2(m+1)}{3} \rfloor} \times2^{\lfloor\frac{3(m+1)-\sigma}{2} \rfloor}$. This contradicts to ${\rm per}D\geq 6^{\lfloor\frac{\sigma-2(m+1)}{3} \rfloor} \times2^{\lfloor\frac{3(m+1)-\sigma}{2} \rfloor}$ in (\ref{equ18}).
Set that $s>1$ or $t>1$.  By (v) of Claim \ref{cla2}, we have ${\rm per}D<6^{\lfloor\frac{\sigma-2(m+1)}{3} \rfloor} \times2^{\lfloor\frac{3(m+1)-\sigma}{2} \rfloor}$, a contradiction.

Assume $k_{2}\neq k_{3}$. We will discuss 5 cases, i.e., $\sigma-2(m+1)\equiv2(mod~3)$, and $3(m+1)-\sigma$ is odd; $\sigma-2(m+1)\equiv2(mod~3)$, and $3(m+1)-\sigma$ is even; $\sigma-2(m+1)\equiv1(mod~3)$, and $3(m+1)-\sigma$ is odd; $\sigma-2(m+1)\equiv1(mod~3)$, and $3(m+1)-\sigma$ is even; $\sigma-2(m+1)\equiv0(mod~3)$.

Set $\sigma-2(m+1)\equiv2(mod~3)$, and $3(m+1)-\sigma$ is odd. If $k_{2},k_{3}\in\{1,2\}$, $s=1$ and $t=2$, then ${\rm per}A={\rm per}D=k_{2}\times6^{\lfloor\frac{\sigma-2(m+1)-s}{3} \rfloor} \times2^{\lfloor\frac{3(m+1)-\sigma+s-1}{2} \rfloor}+k_{3}\times6^{\lfloor\frac{\sigma-2(m+1)-t}{3} \rfloor} \times2^{\lfloor\frac{3(m+1)-\sigma+t-1}{2} \rfloor}=2\times6^{\lfloor\frac{\sigma-2(m+1)}{3} \rfloor} \times2^{\lfloor\frac{3(m+1)-\sigma}{2} \rfloor}<4\times6^{\lfloor\frac{\sigma-2(m+1)}{3} \rfloor} \times2^{\lfloor\frac{3(m+1)-\sigma}{2} \rfloor}$. This contradicts to ${\rm per}D\geq 4\times6^{\lfloor\frac{\sigma-2(m+1)}{3} \rfloor} \times2^{\lfloor\frac{3(m+1)-\sigma}{2} \rfloor}$ in (\ref{equ18}).
 By (i) of Claim \ref{cla2}, we have ${\rm per}D<4\times6^{\lfloor\frac{\sigma-2(m+1)}{3} \rfloor} \times2^{\lfloor\frac{3(m+1)-\sigma}{2} \rfloor}$, a contradiction. If $s\neq1$ or $t\neq2$, then this contradicts to ${\rm per}D\geq 4\times6^{\lfloor\frac{\sigma-2(m+1)}{3} \rfloor} \times2^{\lfloor\frac{3(m+1)-\sigma}{2} \rfloor}$ in (\ref{equ18}) by (i) of Claim \ref{cla2}. If $k_{2}\notin\{1,2\}$ or $k_{3}\notin\{1,2\}$,  then we have ${\rm per}D<4\times6^{\lfloor\frac{\sigma-2(m+1)}{3} \rfloor} \times2^{\lfloor\frac{3(m+1)-\sigma}{2} \rfloor}$ by (\ref{equ18}), a contradiction.

Set $\sigma-2(m+1)\equiv2(mod~3)$, and $3(m+1)-\sigma$ is even. If $k_{2}=\frac{3}{2}$ or $k_{3}=\frac{3}{2}$, then this contradicts to ${\rm per}D\geq 3\times6^{\lfloor\frac{\sigma-2(m+1)}{3} \rfloor} \times2^{\lfloor\frac{3(m+1)-\sigma}{2} \rfloor}$ in (\ref{equ18}). If $k_{2}\neq\frac{3}{2}$ or $k_{3}\neq\frac{3}{2}$, then this contradicts to ${\rm per}D\geq 3\times6^{\lfloor\frac{\sigma-2(m+1)}{3} \rfloor} \times2^{\lfloor\frac{3(m+1)-\sigma}{2} \rfloor}$ in (\ref{equ18}) by (v) of Claim \ref{cla2}.

Set $\sigma-2(m+1)\equiv1(mod~3)$, and $3(m+1)-\sigma$ is odd. If $k_{2},k_{3}\in\{1,3\}$, $s=2$ and $t=1$, then this contradicts to ${\rm per}D\geq 2\times6^{\lfloor\frac{\sigma-2(m+1)}{3} \rfloor} \times2^{\lfloor\frac{3(m+1)-\sigma}{2} \rfloor}$ in (\ref{equ18}).  If $s\neq1$ or $t\neq2$, by (iv) of Claim \ref{cla2}, this contradicts to ${\rm per}D\geq 2\times6^{\lfloor\frac{\sigma-2(m+1)}{3} \rfloor} \times2^{\lfloor\frac{3(m+1)-\sigma}{2} \rfloor}$ in (\ref{equ18}).  If $k_{2}\notin\{1,3\}$ or $k_{3}\notin\{1,3\}$, then this contradicts to ${\rm per}D\geq 2\times6^{\lfloor\frac{\sigma-2(m+1)}{3} \rfloor}$ again by (v) of Claim \ref{cla2}.

Set $\sigma-2(m+1)\equiv1(mod~3)$, and $3(m+1)-\sigma$ is even.  If $k_{1},k_{2},k_{3}$ do not satisfy the conditions of the cases of \{$k_{2}=k_{3}=2$ and $\sigma-2(m+1)\equiv2(mod~3)$, and $3(m+1)-\sigma$ is odd; $k_{2}=k_{3}=1$ and  $\sigma-2(m+1)\equiv2(mod~3)$, and $3(m+1)-\sigma$ is odd; $k_{2}=k_{3}=\frac{3}{2}$  and $\sigma-2(m+1)\equiv1(mod~3)$, and $3(m+1)-\sigma$ is odd; $k_{2}=k_{3}=3$ and $\sigma-2(m+1)\equiv1(mod~3)$, and $3(m+1)-\sigma$ is odd; $k_{2}=k_{3}=1$ and $\sigma-2(m+1)\equiv1(mod~3)$, and $3(m+1)-\sigma$ is odd; $k_{2}=k_{3}=1$ and $\sigma-2(m+1)\equiv1(mod~3)$, and $3(m+1)-\sigma$ is even; $k_{2}=k_{3}=4$ and $\sigma-2(m+1)\equiv0(mod~3)$\} in $(m+1)$-square  matrix $D$ as above, then $k_{2}\times6^{\lfloor\frac{\sigma-2(m+1)-s}{3} \rfloor} \times2^{\lfloor\frac{3(m+1)-\sigma+s-1}{2} \rfloor}<\frac{1}{2}\times \frac{3}{2}\times6^{\lfloor\frac{\sigma-2(m+1)-s}{3} \rfloor} \times2^{\lfloor\frac{3(m+1)-\sigma+s-1}{2} \rfloor}$. Similarly, we have $k_{3}\times6^{\lfloor\frac{\sigma-2(m+1)-s}{3} \rfloor} \times2^{\lfloor\frac{3(m+1)-\sigma+s-1}{2} \rfloor}<\frac{1}{2}\times \frac{3}{2}\cdot6^{\lfloor\frac{\sigma-2(m+1)-s}{3} \rfloor} \times2^{\lfloor\frac{3(m+1)-\sigma+s-1}{2} \rfloor}$. Since $k_{2}\neq k_{3}$,  we get that $k_{2}\times6^{\lfloor\frac{\sigma-2(m+1)-s}{3} \rfloor} \times2^{\lfloor\frac{3(m+1)-\sigma+s-1}{2} \rfloor}+k_{3}\times6^{\lfloor\frac{\sigma-2(m+1)-t}{3} \rfloor} \times2^{\lfloor\frac{3(m+1)-\sigma+t-1}{2} \rfloor}<\frac{1}{2}\times \frac{3}{2}\times6^{\lfloor\frac{\sigma-2(m+1)-s}{3} \rfloor} \times2^{\lfloor\frac{3(m+1)-\sigma+s-1}{2} \rfloor}+\frac{1}{2}\times \frac{3}{2}\times6^{\lfloor\frac{\sigma-2(m+1)-s}{3} \rfloor} \times2^{\lfloor\frac{3(m+1)-\sigma+s-1}{2} \rfloor}=\frac{3}{2}\times6^{\lfloor\frac{\sigma-2(m+1)}{3} \rfloor} \times2^{\lfloor\frac{3(m+1)-\sigma}{2} \rfloor}$. This contradicts to ${\rm per}D\geq \frac{3}{2}\times6^{\lfloor\frac{\sigma-2(m+1)}{3} \rfloor} \times2^{\lfloor\frac{3(m+1)-\sigma}{2} \rfloor}$ in (\ref{equ18}).

Set $\sigma-2(m+1)\equiv0(mod~3)$. If $k_{2}=\frac{3}{2}$ or $k_{3}=\frac{3}{2}$, then this contradicts to ${\rm per}D\geq 6^{\lfloor\frac{\sigma-2(m+1)}{3} \rfloor} \times2^{\lfloor\frac{3(m+1)-\sigma}{2} \rfloor}$ in (\ref{equ18}). If $k_{2}\neq\frac{3}{2}$ or $k_{3}\neq\frac{3}{2}$, then this contradicts to ${\rm per}D\geq 6^{\lfloor\frac{\sigma-2(m+1)}{3} \rfloor} \times2^{\lfloor\frac{3(m+1)-\sigma}{2} \rfloor}$ in (\ref{equ18}) by (v) of Claim \ref{cla2}.

\textbf{Case 2}. Suppose that $\sigma=3m+2$.  Checking $(m+1)$-square  matrix $A$, we know that some line of $(m+1)$-square  matrix $A$ has a single $1$, two $1$'s or three $1$'s.
We consider three subcases as follows.

\textbf{Subcase 2.1}. Assume that some line of $(m+1)$-square  matrix $A$ has exactly two $1$'s. Then the matrix form of $(m+1)$-square  matrix $A$ is equals to  $D=\left[
\begin{matrix}
\begin{array}{c|c|c}

1  &  1  &  0\cdots0 \\
\hline
D_{1}  &  D_{2}  &    C     \\
\end{array}
\end{matrix}\right]$.
The matrix forms of $(m+1)$-square  matrix $D$ is determined, see Case 1.

\textbf{Subcase 2.2}.
Suppose that some line of $(m+1)$-square  matrix $A$ has exactly three $1$'s. Then the matrix form of $(m+1)$-square  matrix $A$ is equals to
$E=\left[
\begin{matrix}
\begin{array}{c|c|c|c}

1  &  1  &  1 &  0\cdots0 \\
\hline
E_{1}  &  E_{2}  &  E_{3} &  C   \\
\end{array}
\end{matrix}\right]$.
Assume that $E_{1}$ has $s$~$1$'s,  $E_{2}$ has $t$~$1$'s and $E_{3}$ has $z$~$1$'s. We consider $5$ subsubcases as follows.

\textbf{Subsubcase 2.2.1}. If $s=0$, $t=0$ or $z=0$. Then the matrix $(m+1)$-square  matrix $E$ must be changed  the matrix $(m+1)$-square  matrix $B$ in Case 1. So, the conclusion holds.

\textbf{Subsubcase 2.2.2}. If $s=1$, $t=1$ or $z=1$. Then the matrix $(m+1)$-square  matrix $E$ must be changed  the $(m+1)$-square  matrix $D$ in Case 1. So, the conclusion holds.

\textbf{Subsubcase 2.2.3}. If $s=2$, $t=2$ or $z=2$. Then  some line of $m$-square matrix $E(1,1)$, $E(1,2)$ or $E(1,3)$ has  two $1$'s. Checking  forms of $m$-square submatrices $E(1,1)$, $E(1,2)$ and $E(1,3)$, we know that $m$-square submatrices $E(1,1)$, $E(1,2)$ or $E(1,3)$  must be changed  the $m$-square matrix $D$ in Case 1. By induction hypothesis,  we can know that $m$-square submatrix $E(1,1)$ is combinatorially equivalent to one of matrices: $N_{m,\sigma}$, $M_{m,\sigma}$, $T_{m,\sigma}$, $P_{m,\sigma}$. Hence, we obtain that $(m+1)$-square matrix $A(=E)$ is combinatorially equivalent to $N_{(m+1),\sigma}$, $M_{(m+1),\sigma}$,  $T_{(m+1),\sigma}$, $P_{(m+1),\sigma}^{*}$, respectively.

\textbf{Subsubcase 2.2.4}. If $s=3$, $t=3$ and $z=3$. Then the number of $1$'s in $m$-square submatrix $E(1,1)$ is $3m-3$.  By induction hypothesis,  we know that ${\rm per}E(1,1)={\rm per}E(1,2)={\rm per}E(1,3)=k_{2}\times6^{\lfloor\frac{3m-3-2m}{3} \rfloor} \times2^{\lfloor\frac{3m-(3m-3)}{2} \rfloor}=\frac{1}{3}\times k_{2}\times6^{\lfloor\frac{m} {3}\rfloor}$.  By Laplace Expansion Theorem,  we have ${\rm per}E=3\times {\rm per}E(1,1)= k_{2}\times6^{\lfloor\frac{m} {3}\rfloor}$. We further consider the following three cases. Assume that $(m+1)\equiv0(mod~3)$. The permanent of each  matrix of the form $F_{(m+1),\sigma}$, $U_{(m+1),\sigma}$, $V_{(m+1),\sigma}$, $W_{(m+1),\sigma}$, $X_{(m+1),\sigma}$, $Y_{(m+1),\sigma}$ and $F_{(m+1),\sigma}^{*}$ equals to $k_{1}\times6^{\lfloor\frac{(m+1)}{3} \rfloor}$, where $k_{1}=1$. By $(m+1)\equiv0(mod~3)$, we obtain that $m\equiv2(mod~3)$.  This means that $k_{2}\in\{3,4\}$ in ${\rm per}E(1,1)$. Hence, ${\rm per}E= k_{2}\times6^{\lfloor\frac{(m+1)-1} {3}\rfloor}=k_{2}\times6^{\lfloor\frac{m+1} {3}\rfloor}\times\frac{1}{6}\leq \frac{2}{3}\times6^{\lfloor\frac{(m+1)} {3}\rfloor}<1\times6^{\lfloor\frac{(m+1)} {3}\rfloor}$, which contradicts to the fact that  ${\rm per}A=\mu((m+1),\tau)$. Assume that $(m+1)\equiv1(mod~3)$. The permanent of each matrix of the form $R_{(m+1),\sigma}$, $S_{(m+1),\sigma}$, $T_{(m+1),\sigma}$, $Q_{(m+1),\sigma}$, $P_{(m+1),\sigma}^{*}$ is $k_{1}\times6^{\lfloor\frac{(m+1)}{3} \rfloor}$, where $k_{1}\in\{\frac{3}{2},2\}$. By $(m+1)\equiv1(mod~3)$, we get that  $m\equiv0(mod~3)$. This implies that $k_{2}=1$ in ${\rm per}E(1,1)$. Hence, ${\rm per}E= k_{2}\times6^{\lfloor\frac{(m+1)-1} {3}\rfloor}=1\times6^{\lfloor\frac{(m+1)} {3}\rfloor}\times\frac{1}{6}< \frac{3}{2}\times6^{\lfloor\frac{(m+1)} {3}\rfloor}$, which contradicts to the fact that  ${\rm per}E=\mu((m+1),\tau)$. Assume that $(m+1)\equiv2(mod~3)$. The permanent of each matrix of the form $M_{(m+1),\sigma}$, $N_{(m+1),\sigma}$ is $k_{1}\times6^{\lfloor\frac{(m+1)}{3} \rfloor}$, where $k_{1}\in\{3,4\}$. By $(m+1)\equiv2(mod~3)$, we have $m\equiv1(mod~3)$. This means that $k_{2}\in\{\frac{3}{2},2\}$. Hence, ${\rm per}E= k_{2}\times6^{\lfloor\frac{(m+1)-1} {3}\rfloor}\leq2\times6^{\lfloor\frac{(m+1)} {3}\rfloor}< 3\times6^{\lfloor\frac{(m+1)} {3}\rfloor}$. This contradicts to ${\rm per}E\geq k_{1}\times6^{\lfloor\frac{(m+1)}{3} \rfloor}$.

\textbf{Subsubcase 2.2.5}. If  $s>3$, $t>3$ and $z>3$. Then we can obtain that the number of $1$'s is  $3m-1-s$ in $m$-square submatrix $E(1,1)$.  We further consider the following three cases. Assume that $m\equiv0(mod~3)$.  By induction hypothesis,  set ${\rm per}E(1,1)=k_{2}\times6^{\lfloor\frac{\sigma-2m-s-3}{3} \rfloor} \times2^{\lfloor\frac{3m-\sigma+s+3}{2} \rfloor}=k_{2}\times6^{\lfloor\frac{m-1-s}{3}\rfloor} \times2^{\lfloor\frac{1+s}{2}\rfloor}\leq \frac{1}{3}\times k_{2}\times6^{\lfloor\frac{m}{3} \rfloor} \times6^{\frac{-1-s}{3}} \times2^{\frac{s+1}{2}}=\frac{1}{3}\times k_{2}\times6^{\lfloor\frac{m}{3} \rfloor}\times(\frac{\sqrt{2}}{\sqrt[3]{6}})^{s+1}$, where $k_{2}=1$. Since $\frac{\sqrt{2}}{\sqrt[3]{6}}< 1$, we know that if $s$ is enlarged then $k_{2}\times6^{\lfloor\frac{\sigma-2m-s-3}{3} \rfloor} \times2^{\lfloor\frac{3m-\sigma+s+3}{2} \rfloor}$ is lessened. Hence, we get that ${\rm per}E(1,1)<\frac{1}{3}\times k_{2}\times6^{\lfloor\frac{m}{3}\rfloor}$, ${\rm per}E(1,2)<\frac{1}{3}\times k_{2}\times6^{\lfloor\frac{m} {3}\rfloor}$ and ${\rm per}E(1,3)<\frac{1}{3}\times k_{2}\times6^{\lfloor\frac{m} {3}\rfloor}$.   By Laplace Expansion Theorem, we have we have ${\rm per}E< k_{2}\times6^{\lfloor\frac{m} {3}\rfloor}$. Then this contradicts to ${\rm per}E\geq k_{1}\times6^{\lfloor\frac{(m+1)}{3} \rfloor}$ by \textbf{Subsubcase 2.2.4} above. Assume that $m\equiv1(mod~3)$.  By induction hypothesis,  set ${\rm per}E(1,1)=k_{2}\times6^{\lfloor\frac{\sigma-2m-s-3}{3} \rfloor} \times2^{\lfloor\frac{3m-\sigma+s+3}{2} \rfloor}=k_{2}\times6^{\lfloor\frac{m-1-s}{3}\rfloor} \times2^{\lfloor\frac{1+s}{2} \rfloor}\leq \frac{1}{3}\times k_{2}\times6^{\lfloor\frac{m-1}{3} \rfloor} \times6^{\frac{-s}{3}} \times2^{\frac{s}{2}}=\frac{1}{3}\times k_{2}\times6^{\lfloor\frac{m-1}{3} \rfloor}\times(\frac{\sqrt{2}}{\sqrt[3]{6}})^{s}$, where $k_{2}\in\{\frac{3}{2},2\}$. Since $\frac{\sqrt{2}}{\sqrt[3]{6}}< 1$, we know that if $s$ is enlarged then $k_{2}\times6^{\lfloor\frac{\sigma-2m-s-3}{3} \rfloor} \times2^{\lfloor\frac{3m-\sigma+s+3}{2} \rfloor}$ is lessened. Hence, we get that ${\rm per}E(1,1)<\frac{1}{3}\times k_{2}\times6^{\lfloor\frac{m-1}{3}\rfloor}$, ${\rm per}E(1,2)<\frac{1}{3}\times k_{2}\times6^{\lfloor\frac{m-1} {3}\rfloor}$ and ${\rm per}E(1,3)<\frac{1}{3}\times k_{2}\times6^{\lfloor\frac{m-1} {3}\rfloor}$.   By Laplace Expansion Theorem, we have we have ${\rm per}E< k_{2}\times6^{\lfloor\frac{m-1} {3}\rfloor}$. Then this contradicts to ${\rm per}E\geq k_{1}\times6^{\lfloor\frac{(m+1)}{3} \rfloor}$ by \textbf{Subsubcase 2.2.4} above. Assume that $m\equiv2(mod~3)$.  By induction hypothesis,  set ${\rm per}E(1,1)=k_{2}\times6^{\lfloor\frac{\sigma-2m-s-3}{3} \rfloor} \times2^{\lfloor\frac{3m-\sigma+s+3}{2} \rfloor}=k_{2}\times6^{\lfloor\frac{m-1-s}{3}\rfloor} \times2^{\lfloor\frac{1+s}{2} \rfloor}\leq \frac{1}{3}\times k_{2}\times6^{\lfloor\frac{m-2}{3} \rfloor} \times6^{\frac{1-s}{3}} \times2^{\frac{s-1}{2}}\times2=\frac{2}{3}\times k_{2}\times6^{\lfloor\frac{m-1}{3} \rfloor}\times(\frac{\sqrt{2}}{\sqrt[3]{6}})^{s-1}$, where $k_{2}\in\{3,4\}$. Since $\frac{\sqrt{2}}{\sqrt[3]{6}}< 1$, we know that if $s$ is enlarged then $k_{2}\times6^{\lfloor\frac{\sigma-2m-s-3}{3} \rfloor} \times2^{\lfloor\frac{3m-\sigma+s+3}{2} \rfloor}$ is lessened. Hence, we get that ${\rm per}E(1,1)<\frac{1}{3}\times k_{2}\times6^{\lfloor\frac{m-2}{3}\rfloor}$, ${\rm per}E(1,2)<\frac{1}{3}\times k_{2}\times6^{\lfloor\frac{m-2} {3}\rfloor}$ and ${\rm per}E(1,3)<\frac{1}{3}\times k_{2}\times6^{\lfloor\frac{m-2} {3}\rfloor}$.   By Laplace Expansion Theorem, we have we have ${\rm per}E< k_{2}\times6^{\lfloor\frac{m-2} {3}\rfloor}$. Then this contradicts to ${\rm per}E\geq k_{1}\times6^{\lfloor\frac{(m+1)}{3} \rfloor}$ by \textbf{Subsubcase 2.2.4} above. Similarly, the conclusion holds when  $m\equiv1(mod~3)$ or $m\equiv2(mod~3)$.

\textbf{Subcase 2.3}.
Suppose that some line of $(m+1)$-square  matrix $A$ has exactly one $1$. Then there exists other line of $(m+1)$-square matrix $A$ has $q(\geq 4)$~$1$'s, and the  form of $(m+1)$-square  matrix $A$ is equals to matrix $H=\left[\begin{matrix}
\begin{array}{c|c|c|c|c}

1  &  1  & \cdots& 1 &  0\cdots0 \\
\hline
H_{1}  &  H_{2} & \cdots &  H_{q} &  C   \\
\end{array}
\end{matrix}\right].$

\textbf{Subsubcase 2.3.1}. If $s=0$, $t=0$ or $z=0$. Then the matrix $(m+1)$-square  matrix $H$ must be changed  the $(m+1)$-square  matrix $B$ in Case 1. So, the conclusion holds.

\textbf{Subsubcase 2.3.2}. If $s=1$, $t=1$ or $z=1$. Then the $(m+1)$-square  matrix $H$ must be changed  the $(m+1)$-square  matrix $D$ in Case 1. So, the conclusion holds.

\textbf{Subsubcase 2.3.3}. If $s=2$, $t=2$ or $z=2$. hen the $(m+1)$-square  matrix $H$ must be changed  the $(m+1)$-square  matrix $E$ in Case 2. So, the conclusion holds.

\textbf{Subsubcase 2.3.4}. If $s\geq3$, $t\geq3$ and $z\geq3$. Then we can obtain that the number of $1$'s is  $3m+2-q-s$ in $m$-square submatrix $H(1,1)$. If $n=m+1$ then the permanents of
$F_{m+1,\sigma}$, $F_{m+1,\sigma}^{*}$, $U_{m+1,\sigma}$, $V_{m+1,\sigma}$, $W_{m+1,\sigma}$, $X_{m+1,\sigma}$, $Y_{m+1,\sigma}$, $Y_{m+1,\sigma}$, $F_{m+1,\sigma}^{*}$, $R_{m+1,\sigma}$, $S_{m+1,\sigma}$, $T_{m+1,\sigma}$,  $Q_{m+1,\sigma}$, $P_{m+1,\sigma}^{*}$, $M_{m+1,\sigma}$ and $N_{m+1,\sigma}$  equal to  $k_{1}\times6^{\lfloor\frac{\sigma-2(m+1)}{3} \rfloor} \times2^{\lfloor\frac{3(m+1)-\sigma}{2} \rfloor}=k_{1}\times6^{\lfloor\frac{m}{3} \rfloor}$,  where $k_{1}\in\{1, \frac{3}{2}, 2, 3,  4\}$. These imply that ${\rm per}H=\mu(m+1,\tau)\geq k_{1}\times6^{\lfloor\frac{\sigma-2(m+1)}{3} \rfloor} \times2^{\lfloor\frac{3(m+1)-\sigma}{2} \rfloor}$, where $k_{1}\in\{1, \frac{3}{2}, 2, 3,  4\}$.  Assume that $m\equiv0(mod~3)$. By induction hypothesis, we can know that ${\rm per}H(1,1)=k_{2}\times6^{\lfloor\frac{m+2-q-s} {3}\rfloor}\times2^{\lfloor\frac{q+s-2}{2} \rfloor}=k_{2}\times6^{\lfloor\frac{m} {3}\rfloor}\times6^{\lfloor\frac{2-q-s} {3}\rfloor}\times2^{\lfloor\frac{q+s-2}{2} \rfloor}\leq k_{2}\times6^{\lfloor\frac{m} {3}\rfloor}\times6^{\frac{2-q-s} {3}}\times2^{\frac{q+s-2}{2}}=k_{2}\times6^{\lfloor\frac{m} {3}\rfloor}\times(\frac{\sqrt{2}}{\sqrt[3]{6}})^{q+s-2}\leq k_{2}\times6^{\lfloor\frac{m} {3}\rfloor}\times(\frac{\sqrt{2}}{\sqrt[3]{6}})^{5}$. Since $\frac{\sqrt{2}}{\sqrt[3]{6}}<1$, we know that if $q$ or $s$ is enlarged then $k_{2}\times6^{\lfloor\frac{m+2-q-s} {3}\rfloor}\times2^{\lfloor\frac{q+s-2}{2} \rfloor}$ is lessened. Hence, ${\rm per}H(1,1)=k_{2}\times6^{\lfloor\frac{m+2-q-s} {3}\rfloor}\times2^{\lfloor\frac{q+s-2}{2} \rfloor}\leq k_{2}\times6^{\lfloor\frac{m+2-4-3} {3}\rfloor}\times2^{\lfloor\frac{4+3-2}{2} \rfloor}=k_{2}\times6^{\lfloor\frac{m}{3} \rfloor}=\frac{1}{9}\times k_{2}\times6^{\lfloor\frac{m}{3} \rfloor}<\frac{1}{3}\times k_{2}\times6^{\lfloor\frac{m} {3}\rfloor}$, where $k_{2}\in\{1, \frac{3}{2}, 2, 3,  4\}$.  Then this contradicts to ${\rm per}H\geq k_{1}\times6^{\lfloor\frac{m}{3} \rfloor}$ by \textbf{Subsubcase 2.2.5} above.  Assume that $m\equiv1(mod~3)$. By induction hypothesis, we can know that ${\rm per}H(1,1)=k_{2}\times6^{\lfloor\frac{m+2-q-s} {3}\rfloor}\times2^{\lfloor\frac{q+s-2}{2} \rfloor}=k_{2}\times6^{\lfloor\frac{m} {3}\rfloor}\times6^{\lfloor\frac{2-q-s} {3}\rfloor}\times2^{\lfloor\frac{q+s-2}{2} \rfloor}\leq k_{2}\times6^{\lfloor\frac{m} {3}\rfloor}\times6^{\frac{2-q-s} {3}}\times2^{\frac{q+s-2}{2}}=k_{2}\times6^{\lfloor\frac{m} {3}\rfloor}\times(\frac{\sqrt{2}}{\sqrt[3]{6}})^{q+s-2}\leq k_{2}\times6^{\lfloor\frac{m} {3}\rfloor}\times(\frac{\sqrt{2}}{\sqrt[3]{6}})^{5}$. Since $\frac{\sqrt{2}}{\sqrt[3]{6}}<1$, we know that if $q$ or $s$ is enlarged then $k_{2}\times6^{\lfloor\frac{m+2-q-s} {3}\rfloor}\times2^{\lfloor\frac{q+s-2}{2} \rfloor}$ is lessened. Hence, ${\rm per}H(1,1)=k_{2}\times6^{\lfloor\frac{m+2-q-s} {3}\rfloor}\times2^{\lfloor\frac{q+s-2}{2} \rfloor}\leq k_{2}\times6^{\lfloor\frac{m+2-4-3} {3}\rfloor}\times2^{\lfloor\frac{4+3-2}{2} \rfloor}=k_{2}\times6^{\lfloor\frac{m}{3} \rfloor}=\frac{1}{9}\times k_{2}\times6^{\lfloor\frac{m}{3} \rfloor}<\frac{1}{3}\times k_{2}\times6^{\lfloor\frac{m} {3}\rfloor}$, where $k_{2}\in\{1, \frac{3}{2}, 2, 3,  4\}$.  Then this contradicts to ${\rm per}H\geq k_{1}\times6^{\lfloor\frac{m}{3} \rfloor}$ by \textbf{Subsubcase 2.2.5} above. Assume that $m\equiv2(mod~3)$. By induction hypothesis, we can know that ${\rm per}H(1,1)=k_{2}\times6^{\lfloor\frac{m+2-q-s} {3}\rfloor}\times2^{\lfloor\frac{q+s-2}{2} \rfloor}=k_{2}\times6^{\lfloor\frac{m} {3}\rfloor}\times6^{\lfloor\frac{2-q-s} {3}\rfloor}\times2^{\lfloor\frac{q+s-2}{2} \rfloor}\leq k_{2}\times6^{\lfloor\frac{m} {3}\rfloor}\times6^{\frac{2-q-s} {3}}\times2^{\frac{q+s-2}{2}}=k_{2}\times6^{\lfloor\frac{m} {3}\rfloor}\times(\frac{\sqrt{2}}{\sqrt[3]{6}})^{q+s-2}\leq k_{2}\times6^{\lfloor\frac{m} {3}\rfloor}\times(\frac{\sqrt{2}}{\sqrt[3]{6}})^{5}$. Since $\frac{\sqrt{2}}{\sqrt[3]{6}}<1$, we know that if $q$ or $s$ is enlarged then $k_{2}\times6^{\lfloor\frac{m+2-q-s} {3}\rfloor}\times2^{\lfloor\frac{q+s-2}{2} \rfloor}$ is lessened. Hence, ${\rm per}H(1,1)=k_{2}\times6^{\lfloor\frac{m+2-q-s} {3}\rfloor}\times2^{\lfloor\frac{q+s-2}{2} \rfloor}\leq k_{2}\times6^{\lfloor\frac{m+2-4-3} {3}\rfloor}\times2^{\lfloor\frac{4+3-2}{2} \rfloor}=k_{2}\times6^{\lfloor\frac{m}{3} \rfloor}=\frac{1}{9}\times k_{2}\times6^{\lfloor\frac{m}{3} \rfloor}<\frac{1}{3}\times k_{2}\times6^{\lfloor\frac{m} {3}\rfloor}$, where $k_{2}\in\{1, \frac{3}{2}, 2, 3,  4\}$.  Then this contradicts to ${\rm per}H\geq k_{1}\times6^{\lfloor\frac{m}{3} \rfloor}$ by \textbf{Subsubcase 2.2.5} above.

\textbf{Case 3}. Suppose that $\sigma=3m+3$.  Checking $(m+1)$-square  matrix $A$, it is easy to see that  some line of $(m+1)$-square  matrix $A$ has one $1$, two $1$'s or three $1$'s. If some line of $(m+1)$-square  matrix $A$ has a single $1$ or two $1$'s, then there exists other line of $(m+1)$-square  matrix $A$ has $q(\geq4)$~$1$'s, and the matrix form of $(m+1)$-square  matrix $A$  equals to matrix
$H=\left[\begin{matrix}
\begin{array}{c|c|c|c|c}

1  &  1  & \cdots& 1 &  0\cdots0 \\
\hline
H_{1}  &  H_{2} & \cdots &  H_{q} &  C   \\
\end{array}
\end{matrix}\right].$
If some line of $(m+1)$-square  matrix $A$ has exectly three $1$'s, then the matrix form of $(m+1)$-square  matrix $A$  equals to matrix
$E=\left[
\begin{matrix}
\begin{array}{c|c|c|c}

1  &  1  &  1 &  0\cdots0 \\
\hline
E_{1}  &  E_{2}  &  E_{3} &  C   \\
\end{array}
\end{matrix}\right]$.
The proof is the same as that of  $\sigma=3m+2$.

Combining the arguments above, the assertion holds.
\qquad \qquad\qquad \quad\quad\qquad \qquad\qquad\qquad$\Box$

\section{ Discussion}

In this paper, we  partially answered  Brualdi-Goldwasser-Michael problem. We characterized the maximum permanents of all matrices in $\mathscr{U}(n,\tau)$  when  [$2n+1$, $3n$]. And  we also characterized the maximum permanents of   all matrices in $\mathscr{U}(n,\tau)$  when $\sigma-kn\equiv0 (mod~k+1)$ and $(k+1)n-\sigma\equiv0(mod~k)$, where $\sigma=n^{2}-\tau$,  $kn\leq\sigma\leq (k+1)n$ and $k$ is integer.
It is very difficult to characterize the maximum permanents of all matrices in $\mathscr{U}(n,\tau)$  when  [$3n+1$,  $n^{2}-2n$]  by using the method similar to this paper. We will look for new ways to characterize the maximum permanents of all matrices in $\mathscr{U}(n,\tau)$  when  [$3n+1$,  $n^{2}-2n$].

\end{document}